\documentclass[11pt]{amsart}

\usepackage{amssymb, tikz, ifthen, url, natbib}

\RequirePackage{ifthen, url}

\newboolean{oneResultSequence}
\setboolean{oneResultSequence}{false} 
\newcommand{\numberResults}{
\ifthenelse{\boolean{oneResultSequence}}{
\newtheorem{axiom}{Axiom}
\newtheorem{theorem}{Theorem}
\newtheorem{lemma}[theorem]{Lemma}
\newtheorem{proposition}[theorem]{Proposition}
\newtheorem{corollary}[theorem]{Corollary}
\newtheorem{clame}[theorem]{Claim}
\theoremstyle{definition}
\newtheorem{definition}[theorem]{Definition}
\newtheorem{example}[theorem]{Example}
\newtheorem{xca}[theorem]{Exercise}
\theoremstyle{remark}
\newtheorem{remark}[theorem]{Remark}
\newtheorem{conjecture}[theorem]{Conjecture}
}{

\newtheorem{theorem}{Theorem}
\newtheorem{lemma}{Lemma}
\newtheorem{proposition}{Proposition}
\newtheorem{corollary}{Corollary}
\newtheorem{clame}{Claim}
\theoremstyle{definition}

\newtheorem{example}{Example}

\theoremstyle{remark}

\newtheorem{conjecture}{Conjecture}
}}

\newcommand{\rslt}[1]{\ifthenelse{\ref{lem:#1} > 0}{lemma \ref{lem:#1}}%
{%
\ifthenelse{\ref{prp:#1} > 0}{proposition \ref{prp:#1}}%
{%
\ifthenelse{\ref{thm:#1} > 0}{theorem \ref{thm:#1}}%
{%
\ifthenelse{\ref{cor:#1} > 0}{corollary \ref{cor:#1}}%
{%
\ifthenelse{\ref{clm:#1} > 0}{claim \ref{clm:#1}}%
{%
\ifthenelse{\ref{cnj:#1} > 0}{conjecture \ref{cnj:#1}}%
{%
\ifthenelse{\ref{xmpl:#1} > 0}{example \ref{xmpl:#1}}{\textbf{result
    ??} \typeout{Result #1 is not defined.}}%
}%
}%
}%
}%
}%
}%
}

\newcommand{\Rslt}[1]{\ifthenelse{\ref{lem:#1} > 0}{Lemma \ref{lem:#1}}%
{%
\ifthenelse{\ref{prp:#1} > 0}{Proposition \ref{prp:#1}}%
{%
\ifthenelse{\ref{thm:#1} > 0}{Theorem \ref{thm:#1}}%
{%
\ifthenelse{\ref{cor:#1} > 0}{Corollary \ref{cor:#1}}
{%
\ifthenelse{\ref{clm:#1} > 0}{Claim \ref{clm:#1}}%
{%
\ifthenelse{\ref{cnj:#1} > 0}{Conjecture \ref{cnj:#1}}%
{%
\ifthenelse{\ref{xmpl:#1} > 0}{Example \ref{xmpl:#1}}{\textbf{Result
    ??} \typeout{Result #1 is not defined.}}%
}%
}%
}%
}%
}%
}%
}

\newcommand{\eqnref}[1]{\ref{eqn:#1}}
 
\newcommand{\secref}[1]{\ref{sec:#1}}

\newcommand{\lema}[2]{\begin{lemma} #2 \label{lem:#1} \end{lemma}}
\newcommand{\lem}[1]{{l}emma \ref{lem:#1}}

\newcommand{\claim}[2]{\begin{clame} #2 \label{clm:#1} \end{clame}}

\newcommand{\corol}[2]{\begin{corollary} #2 \label{cor:#1} \end{corollary}}
\newcommand{\cor}[1]{{c}orollary \ref{cor:#1}}

\newcommand{\propo}[2]{\begin{proposition} #2 \label{prp:#1} \end{proposition}}
\newcommand{\prp}[1]{{p}roposition \ref{prp:#1}}

\newcommand{\exmpl}[2]{\begin{example} #2 \label{xmpl:#1} \end{example}}

\newcommand{\display}[2]{\begin{equation} #2 \label{eqn:#1} \end{equation}}
\newcommand{\eqn}[1]{(\ref{eqn:#1})}

\newcommand{\Section}[2]{\section{#2\label{sec:#1}}}
\newcommand{\Subsection}[2]{\subsection{#2\label{sec:#1}}}
\renewcommand{\sec}[1]{{s}ection \ref{sec:#1}}

\newenvironment{enum} 
{\begin{list}{\makebox[\labelwidth][l]{(\arabic{enumi})}}{\usecounter{enumi}}
\setcounter{enumi}{\value{equation}}}
{\setcounter{equation}{\value{enumi}} \end{list}}

\newcommand{\meti}[2]{\item #2 \label{eqn:#1}} 

\newcommand{\abbrevEnvir}{
\expandafter\newcommand\expandafter{\csname bi\endcsname}{\begin{itemize}} 
\expandafter\newcommand\expandafter{\csname ei\endcsname}{\end{itemize}}
\expandafter\newcommand\expandafter{\csname be\endcsname}{\begin{enumerate}} 
\expandafter\newcommand\expandafter{\csname ee\endcsname}{\end{enumerate}}
\expandafter\newcommand\expandafter{\csname bc\endcsname}{\begin{center}} 
\expandafter\newcommand\expandafter{\csname ec\endcsname}{\end{center}}
}

\newcommand{\from}{\colon}
\newcommand{\Nat}{\mathbb{N}}
\renewcommand{\Re}{\mathbb{R}}

\author[F.~Borhani]{Fatemeh Borhani}
\email{fatemeborhani@gmail.com}

\author[E.~J.~Green]{Edward J. Green}
\address{Department of Economics, Pennsylvania State University}
\email{eug2@psu.edu}

\newcommand{\hcf}{the Human Capital
  Foundation, (\url{http://www.hcfoundation.ru/en/}) and
  particularly Andrey P.\ Vavilov, for research support through the
  Center for the Study of Auctions, Procurements, and Competition
  Policy (CAPCP, \url{http://capcp.psu.edu/}) at the Pennsylvania
  State University.\ }

\newcommand{\q}[2]{q_{#1#2}}

\newcommand{\restrict}{\upharpoonright}

\newcommand{\depth}{d}
\newcommand{\num}{\nu}
\newcommand{\badchoice}[1][\zeta]{{B_{#1}}} 
\newcommand{\sample}[1][\relax]{{\Omega_{#1}}} 
\newcommand{\leaves}{\Lambda}

\newcommand{\mand}{\text{\ and\ }}

\abbrevEnvir

\renewcommand{\Pr}{\pi}

\newcommand{\mean}{\mathsf{E}}

\newcommand{\br}[1][\relax]{\mathfrak{b}_{#1}}
\newcommand{\avoid}[1][\zeta]{{\beta_{#1}}} 

\newcommand\B{\mathcal{B}}

\newcommand\X{\mathbb{X}}

\DeclareMathOperator{\wms}{\sigma} 
\DeclareMathOperator{\sms}{\Sigma} 
\DeclareMathOperator{\immms}{\Sigma_1} 

\newcommand{\nothing}{\mathsf{r}} 
\newcommand{\incompat}{\perp} 

\newcommand{\set}{{\mathfrak l}}

\newcommand{\event}[1][\relax]{{\mathfrak{e}_{#1}}}

\newcommand{\pref}{R}
\newcommand{\immed}[1]{Y(#1)} 

\newcommand{\sceup}{{S}CEU rational} 
\newcommand{\sceu}{{S}CEU-rational\ }

\newcommand{\domain}{B}
\newcommand{\range}{\mathcal{S}}
\newcommand{\graph}{\Omega}
\newcommand{\atom}{\beta}

\providecommand{\fgr}[1]{{f}igure \ref{fig:#1}}

\begin{document}

\numberResults

\title[Evidence-based choice]{A parsimonious theory\\ of
  evidence-based choice} \date{2018.02.21\\ \strut} \thanks{The
  authors gratefully acknowledge \hcf Ed Green thanks CEMFI and its
  faculty and staff for their hospitality at the inception of his
  research.  Fatemeh Borhani's research was conducted while she was a
  graduate student at the Pennsylvania State University and a visiting
  assistant professor at the University of Pittsburgh. Both authors
  appreciate the generous and insightful advice received from a
  multitude of people during the long gestation of the article.}

\begin{abstract}
That an agent's possible evidential states form a Boolean algebra (on
which it is natural to define a probability measure) is an assertion
that ideally should be proved, rather than assumed, in justifying
rational choice as a representation of expected utility. A more
parsimonious, axiomatic characterization of evidence is provided
here. Two primitive entities are \emph{evidential states\/} and a
relation, \emph{more specific than,} between evidential states. The
axioms specify that more-specific-than is a partial order, there is a
minimally specific e-state, and more-specific-than is a separative
order. \emph{Choice alternatives\/} are another primitive entity. A
\emph{plan\/} is an assignment of a choice alternative to each
evidential state. In general, plans satisfying a version of the
sure-thing principle cannot be rationalized by expected utility. But
there is such a rationalization if the evidential structure is a tree.
\end{abstract}

\maketitle

\Section{1}{Introduction}

A main ``selling point'' of the Bayesian approach to statistics is
that it expands the focus of the subject from a rather narrow study of
numerical data to a broader and more conceptually satisfactory study
of evidence in general. Especially, in the decision-theoretic
formulation of \citet{Savage-1972}, it is a study of the rational use
of evidence to inform practical decisions.\footnote{Savage's main
  innovation was arguably to recast statistics as a guide to making
  \emph{group} decisions. First-edition dates of Savage's book and of
  other landmark references are provided in the bibliography.}  From
that perspective, Savage's modeling of evidence as a Boolean algebra
of events, each of which is a set of possible states of the world, is
a liability. There are intuitive examples (one of which will be
presented below, in \sec{b}) that cast doubt on the idea that the sets
corresponding to evidential states would form a Boolean algebra.  In
view of such examples, the thesis that evidential states are sets
that form precisely the type of a structure on which probability
measures are defined is something that it would be more satisfactory to
prove, rather than to hypothesize, in demonstrating that rational
choice can be represented by expected-utility maximization.

This article addresses the following two questions: (a) How should
evidence be modeled for the purpose of normative theorizing about
rational choice? (b) Given how Bayesian statisticians, and Savage in
particular, have modeled evidence, do their conclusions genuinely
characterize how the rational use of evidence should inform choices or
decisions? The first question does not presume Bayesian theory as a
background or benchmark. The second question is distinct from various
questions about, and criticisms of, Bayesian decision theory that have
been raised antecedently.\footnote{In studying whether or not one or
  another principle resembling Savage's sure-thing principle and
  formulable within the language of the theory to be presented here
  entails a conditional-expected-utility representation of rational
  choice, the authors' purpose is not to privilege that
  representation. In fact, the various generalizations of
  expected-utility theory that have been proposed as alternatives to
  Savage's theory (surveyed, for example, by \citet{Gilboa-2009}) all
  incorporate, with at most modest modifications, Savage's
  formalization of evidence. Parallel to the analysis provided here,
  one could study whether the various postulates of those
  generalizations that play analogous roles to the sure-thing
  principle can be re-stated in the present theory in ways that entail
  the representations of rational choice derived respectively from
  them. A plausible conjecture is that the conclusions would mirror
  those found here: that such representations can be derived under
  special assumptions about the structure of evidence, but not from
  parsimonious assumptions that express the intrinsic structure of
  evidence in general.}

Question (a) will be addressed in \sec{b}, where a parsimonious
axiomatization of evidence is proposed. The axioms are
stated in terms of two primitive entities: a set of a person's
possible \emph{evidential states\/} (or \emph{e-states,} for short)
and a binary relation, \emph{weakly more specific than,} between
e-states. An \emph{evidential structure} comprises both of these
entities. The axioms that characterize an evidential structure are
that more-specific-than is a partial order; that there is a
prior-evidence e-state, than which every other e-state is stritly more
specific; and that more-specific-than is a \emph{separative\/}
order. This last condition means that, if some e-state, $y$, is not
more specific than another, $x$, then there is a third e-state, $z$,
which is strictly more specific than $y$ and is also incompatible with
$x$.\footnote{For example, since having observed
  a flash is not more specific than having observed a red flash, then
  there must be another e-state---having observed a blue flash---that
  is also more specific than observing a flash but that is
  incompatible with having observed a red flash. \label{fn:sep}}

Question (b) will be addressed in sections \secref{z}--\secref{B},
where evidence-based choice will be studied. An evidential structure
and a set of \emph{choice alternatives\/} (\emph{acts}, in Savage's
terminology) will be primitive entities. A \emph{plan\/} will be
defined to be an assignment of a choice alternative (to be interpreted
as being the recommended/chosen alternative) to each evidential
state. A plan thus represents an agent's pattern of evidence-based
choices. What it means for a plan to be justified by subjective
conditional expected utility will be defined. Three candidates for an
axiom characterizing rational choice, of increasing strength, and in
the spirit of Savage's sure-thing principle, will be specified, with
the middle-strength axiom being the most closely analogous to a
condition framed by Savage \citeyearpar[theorem 2.3, section
  2.7]{Savage-1972} that he regards to be a formal expression of the
sure-thing principle.\footnote{The weakest of these axioms is
  reminiscent of conditions of \emph{dynamic consistency} and
  \emph{consequentialism} that have been studied by
  \citet{Weller-1978}, \citet{Hammond-1988}, \citet{Ghirardato-2002},
  and others. In various frameworks in which evidential states are
  taken to be sets of possible worlds that form a Boolean algebra,
  those authors have proved results that relate consistency to a
  subjective-conditional-expected-utility criterion. However, there
  are considerable differences in both framework and proof strategy
  between those research contributions and the present one. The
  results proved here do not seem to be derivable as corollaries of
  antecedent results. \label{fn:dynamic}} It will be shown by
counter-example (\rslt G) that not even the strongest of these axioms,
in the presence of the axioms defining an evidential structure,
entails a subjective-conditional-expected-utility rationalization of a
plan.\footnote{By rationalization, it is meant only that some imputation
  of subjective probabilities and state-contingent utilities of acts
  to the agent rationalizes the plan. This is a weaker criterion than
  the criterion of representability that is common in the
  Bayesian-decision-theory literature, which includes an assertion of
  uniqueness (up to positive affine transformation of
  utilities). Antecedent theories of subjective conditional expected
  utility are discussed in \sec{g}.} It will also be proved (\rslt{H}
and \rslt{Q}) that even the weakest version of a sure-thing axiom is
sufficient for a plan to have a
subjective-conditional-expected-utility representation if, besides the
defining axioms of an evidential structure being satisfied, the
specificity relation makes the set of e-states to be a directed tree
(to be defined formally in \eqn{50}). It will be clear that such
evidential structures correspond to filtrations, which play a
prominent role in stochastic-process theory and in applied decision
theory. However, the same intuitive example that shows why it is
restrictive to posit that evidential states form a Boolean algebra,
also shows that the evidential states may not form a tree.

The concluding \sec{c} includes remarks on antecedent theories of
subjective conditional expected utility, the robustness of the
article's results to considering finitely additive measures instead of
countably additive ones, an indication of how weak the parsimonious
theory is, the relationship between infinite evidential structures and
Boolean algebras, and the status of the definition of subjective conditional
expected utility in the article as a model of na\"{i}ve, rather than
sophisticated, Bayesian rationality.

The authors are inclined to give a nuanced interpretation to
the overall results of this article. If the fundamental intuition
about rational choice under uncertainty is that a rational agent
conforms to some version of the sure-thing principle, then subjective
conditional expected utility falls short of being a fully general
necessary condition for rational evidence-based choice, so it should
certainly not be taken to be an intrinsic definition of rational
choice. Nonetheless, its employment for practical decision making is
sound in a fairly general class of situations that can be explicitly,
operationally characterized.

\Section b{Evidence}

\Subsection u{Definition and example of an evidential structure}

The set of an agent's possible \emph{evidential states} will be denoted by
$\X$. (The elements of $\X$ will often be called \emph{e-states,}
for short.) There is a binary relation, $\wms$, among e-states, that
satisfies  assumptions \eqn{5}--\eqn{7} below. That $x \wms y$ ($x$ is
\emph{weakly more specific} than $y$) means that evidence $x$ is consistent with,
and possibly more specific than, evidence $y$. $\sms$ denotes the
strict part of $\wms$ and is defined, as usual, by
\display{1}{x \sms y \iff [x \wms y \text{\ and not\ } y \wms x]}
%
%

E-state $x$ is defined to be immediately more specific
than $z$ ($x \immms z$) by the condition that
\display{2}{x \immms z \iff [x \sms z \text{\ and not\ } \exists y \; x
    \sms y \sms z]} 
Define
\display{3}{\immed z = \{ x \mid x \immms z \}}

Condition \eqn{7}, below, involves the
concept of two e-states being incompatible with one
another. Incompatibility means that there is no e-state that is demonstrably
consistent with both of them, where consistency is demonstrated by
being weakly more specific according to $\wms$. Formally, that $x$
and $y$ are incompatible ($x \incompat y$) is defined by
\display{4}{x \incompat y \iff \text{not\ } \exists z \; [z \wms x
    \mand z \wms y]}

Here are the assumptions regarding e-states and the specificity
ordering. A structure, $(\X, \wms)$, that satisfies conditions
\eqn{5}--\eqn{7} will be called an \emph{evidential structure}, or an
e-structure.
\begin{enum}
\meti{5}{$\wms$ is a partial order (that is, reflexive, transitive,
  and antisymmetric).\footnote{The assumption of antisymmetry could be
    dropped. Then the definition of immediate-dominance
    consistency (cf. \sec l) would have to be amended to require that,
    if $x \wms y$ and $y \wms x$, then the chosen alternative at $x$
    is the same as at $y$. Also, the proofs of several results would
    need to take account of the possibility that $x \wms y$ and $y
    \wms x$ but $x \neq y$. While there can be interpretations of the
    concept of evidence that are intuitively reasonable in some
    applications, and such that the weakly-more-specific relation is
    reflexive and transitive but not antisymmetric (such as the
    example of evidence for identification of a chemical element
    provided in \citet[section 2.1]{bg-2016}), we eschew the highest
    level of generality in order to keep the formal analysis as simple
    as possible.}}
\meti{6} {$\nothing$, an element of $\X$,
  satisfies $\forall x\; [x \neq \nothing \implies x \! \sms
    \nothing]$, and $\X \neq \{ \nothing \}$. (That is, $\nothing$ is
  the root of $(\X, \wms)$, regarded as an ordered structure.)}
\meti{7}{If not $x \wms z$, then there exists an e-state, $y$, such
  that $y \wms x$ and $y \incompat z$.}
\end{enum}

A relation satisfying condition \eqn{7} is called \emph{separative.}
The condition asserts that, unless e-state $x$ is weakly more specific
than e-state $y$, $x$ can be ``fleshed out'' in some way that would be
incompatible with $y$, as has been exemplified in footnote
\ref{fn:sep}.

Although we are not aware of a prior attempt to give an abstract
definition of an evidential structure using assumptions that are as
parsimonious as possible, there are two types of ordered structure
that have been used implicitly for that purpose: Boolean lattices and
trees.\footnote{A Boolean lattice is derived from a Boolean algebra by
  defining $x \wms y \iff x \wedge y = x$. Thus it provides an
  abstract version of the Boolean structure that \citet{Savage-1972}
  imposes on evidence. A tree corresponds to an e-structure that
  satisfies the condition that, for every $x \neq \nothing$, there
  exists a unique $y$ such that $x \immms y$ and there are finitely
  many $y$ such that $x \wms y$. (An equivalent definition is given in
  \eqn{50} below.) \citet{Raiffa-1970} used trees to provide an
  elementary exposition of Bayesian decision theory.} Both of those
types of structure are separative.

An example of an e-structure that is neither a Boolean lattice nor a
tree is now provided, along with two possible interpretations of this
structure. One interpretation (regarding coin tossing) is intended to
provide clear intuition for the meaning of the formal example. The
other interpretation (regarding a scientific investigation) is
intended to convey the breadth of the range of situations in which an
e-structure such as this one might arise. It is not necessarily a
typical e-structure, but it is one that an adequate theory ought to
cover. Thus, the example serves to motivate the study of a more
parsimonious theory of e-structures than the theories of Boolean
lattices and of trees would afford.

\uppercase\fgr a depicts \rslt j, a simple example of a finite
e-structure. In the figure, if one e-state is depicted below another
and the two are connected by a path of line segments, then the lower
e-state is strictly more specific than the higher one. This convention
will be observed in diagrams throughout this article, except where
explicitly amended.

\exmpl j{Let $\X = \{
\nothing, x_1, x_2, z_1, z_2, z_3 \}$. Let $\immms$ comprise the pairs
$(z_1,x_1)$, $(z_2,x_2)$, $(z_3,x_1)$, $(z_3,x_2)$, $(x_1,\nothing)$,
and $(x_2,\nothing)$, and let $\wms$ be the reflexive, transitive
closure of $\immms$. }

\newsavebox{\labelbox}
\sbox{\labelbox}{Diagram of \Rslt j}

\begin{figure} \centering 
\begin{tikzpicture}

\draw (3,3) node[above]{$\nothing = [0,0]$} -- (2,2) node[above left]{$x_1 = [1,0]$};

\draw (3,3) -- (4,2) node[above right]{$x_2 = [0,1]$};

\draw (2,2) -- (1,1) node[below left]{$z_1 = [2,0]$};

\draw (2,2) -- (3,1) node[below]{$z_3 = [1,1]$};

\draw (4,2) -- (5,1) node[below right]{$z_2 = [0,2]$};

\draw (4,2) -- (3,1);

\end{tikzpicture}
\caption{\usebox{\labelbox}} \label{fig:a} \end{figure}

How might this e-structure be interpreted? Suppose that two coins are
tossed, and that the agent is told only a minimum number of tosses
that have landed on each side. That is, such censored reports of the
coin-toss outcomes, rather than the outcomes themselves, correspond to
the agent's possible evidential states. Specifically, let the e-state
associated with $[m,n]$ characterize the agent's situation after
having been told that ``At least $m$ of the coins landed
`heads' and at least $n$ landed `tails'.''

The e-structure in \rslt j is neither a Boolean lattice nor a
tree.\footnote{$(\X, \wms)$ is not a tree because there is not a
  \emph{unique} $y$ such that $z_3 \immms y$. $(\X, \wms)$ cannot be a
  Boolean lattice because, in order for it to be such, there would
  have to be an e-state that is $z_1 \wedge z_2$, that e-state would
  thus be more specific than both $z_1$ and $z_2$, and there is no
  such e-state in $(\X, \wms)$.} Nonetheless, \rslt j is not
pathological. Consider the following analogue. In a certain European
valley that paleontologists have recently begun to explore, Homo
heidelbergensis and Homo neanderthalensis co-existed 300,000 years
ago. A paleontologist wants to know which of the two species was more
numerous at that place and time. The evidential state of the
paleontologist might be, for example, that, so far, 21
H.\ neanderthalensis skeletons and 34 H.\ heidelbergensis skeletons
have been discovered in the valley and shown by radiocarbon dating to
be 300,000 years old. That does not mean that exactly 55 hominins
lived in the valley at the time, of whom 21 were
H.\ neanderthalensis. Rather, it means that there were \emph{at least}
21 members of H.\ neanderthalensis and \emph{at least} 34 members of
H.\ heidelbergensis who lived in the valley then. This is analogous to
the evidential state that at least one toss of a coin landed `heads'
and at least zero tosses landed `tails'. The paleontologist's
evidential state is immediately more specific than each of two
distinct e-states: that 20 H.\ neanderthalensis skeletons and 34
H.\ heidelbergensis skeletons of relevant age have been discovered so
far, and that 21 H.\ neanderthalensis skeletons and 33
H.\ heidelbergensis skeletons of relevant age have been discovered so
far. Thus the e-structure that comprises the possible evidential
states of the paleontologist is not a tree. Neither is
a Boolean lattice. It would be artificial, and arguably
inappropriate, to suppose, for example, that the disjunction of the
paleontologist's e-states is necessarily an e-state.\footnote{In
  particular, the \emph{principle of total evidence} in Bayesian
  decision theory requires that, having made one observation or
  another, an agent should condition on that observation alone rather
  than on its disjunction with an alternate
  observation. Cf.\ \citet{Good-1967}. \label{fn:total}}

\Subsection C{Embedding a finite e-structure in a Boolean field of sets}

As often happens, it is true here that studying the finite examples of
a class of structures can be helpful. All of the examples to be
presented in sections \ref{sec:n} and \ref{sec:y}, which exhibit
disanalogies between a plausible criterion of evidence-based choice
and Savage's \citeyearpar{Savage-1972} Bayesian criterion (that is,
the ``sure-thing principle''), are finite. \Rslt Q (\sec{E}), which
characterizes a case in which a close analogy does hold, is more
easily proved (as \Rslt H in \sec A) in the finite case.
Although a finite e-structure is not a Boolean lattice in general,
\rslt C (to be stated and proved below) will asssert that
every finite e-structure $(\X, \wms)$, corresponds in a canonical way to a finite
Boolean algebra. Specifically, each e-state, $x \in \X$, is associated with the
set of leaves of $\wms$ that are weakly more specific than
$x$. A \emph{leaf} is a maximally specific e-state. That
is,
\begin{enum}
\meti{9}{$x$ is a \emph{leaf} of $(\X, \wms)$ iff not $\exists y \; y \sms
  x$.}
\meti{10}{$\leaves = \{ z \mid z \text{\ is a leaf of\ } (\X, \wms) \}$}
\meti{11}{$\set(x) = \{ z \mid z \in \leaves \mand z \wms x \}$}
\end{enum}
This correspondence will be utilized in the proof of \rslt C, which
proceeds via the following lemma.

\lema D{Let $(\X, \wms)$ be a finite e-structure. For each $x \in \X$,
  there is a leaf, $z$, such that $z \wms x$.}
 
\begin{proof}
Otherwise a contradiction would result, as follows. The e-state $x$
itself cannot be a leaf, or else reflexiveness of $\wms$ (condition
\eqnref{5}) would entail that $x$ itself is a leaf such that is
weakly more specific than $x$. So, there must be
an e-state, $y_0$, such that $y_0 \sms x$. If $y_0$ were a leaf, then
$y_0 \in \set(x)$, so there must be an e-state, $y_1$, such that
$y_1 \sms y_0$. By transitivity of $\sms$ (condition \eqnref{5}),
$y_1 \wms x$, so $y_1$ cannot be a leaf or else $y_1 \in
\set(x)$. Continuing in this way, an infinite sequence $y_0,
y_1,\dotsc, y_n,\dots$ is constructed such that, for each $n$,
$y_{n+1} \sms y_n$. By transitivity and asymmetry of $\sms$, all of
the elements of this sequence are distinct from one another. But this
conclusion contradicts the finiteness of $\X$, so it must be that
$\set(x) \neq \emptyset$.
\end{proof}

\propo C{Let $(\X, \wms)$ be a finite e-structure. The mapping,
  $\set$, defined by \eqn{11}, satisfies
\display{12}{\set(\nothing) = \leaves}
and, for all $x$ and $y$,
\display{13}{\set(x) \neq \emptyset}
\display{14}{x \text{\ is not a leaf\ } \implies \set(x) =
  \bigcup_{z \sms x} \set(z)}
\display{15}{x \wms y \iff \set(x) \subseteq \set(y)}
\display{16}{x \incompat y \iff \set(x) \cap \set(y) = \emptyset}
}

\begin{proof}
Condition \eqn{6}, that every e-state is weakly more specific
than the root, implies condition \eqn{12}.
\Rslt D shows that condition \eqn{13} holds, when \eqn{11}
defines $\set$.

To prove \eqn{14}, suppose that $x$ is not a leaf. If $y \in
\set(x)$, then $y \sms x$ and $y \in \set(y)$, so $y \in
\bigcup_{z \sms x} \set(z)$. Conversely, if $y \in \bigcup_{z \sms
  x} \set(z)$, then for some $z \sms x$, $y \in \set(z)$ and
therefore $y \wms z \wms x$, so $y \in \set(x)$.

To prove condition \eqn{15}, first note that $x \wms y \implies
\set(x) \subseteq \set(y)$ follows from \eqn{14}. To prove the
converse implication, suppose that not $x \wms y$. Then, by
separativeness of $\wms$ (condition \eqn{7}), there exists an
e-state, $w$, such that $w \wms x$ and $w \incompat y$. \uppercase\lem
D shows that, for some leaf, $z$, $z \wms w$. By transitivity of
$\wms$, it follows that $z \wms x$ and $z \incompat y$. Therefore $z
\in \set(x) \setminus \set(y)$, so $\set(x) \not \subseteq
\set(y)$. 

Condition \eqn{16} is proved analogously to \eqn{15}, by proving the
contrapositive of the implication in each direction.
\end{proof}

The conceptual significance of \prp C is that the leaves of $\wms$ can
be treated formally as being states of the world, and the power set of those
leaves can be treated as a Boolean field of events, among which some
(the images under $\set$ of e-states) correspond to evidential
states.\footnote{`Field' will be used in this article to refer either
  to a Boolean field or a $\sigma$ field. In either case, a field is a
  subset of the power set of some specified set. It must contain the
  entire space and the null set, and must be closed under
  complementation and under finite (in the case of a Boolean field) or
  countable (in the case of a $\sigma$ field) unions and intersections
  of its elements. The term will be used without qualification when a
  statement applies to both Boolean fields and $\sigma$ fields, and
  also when the type of field intended is clear from context.} The
one-to-one mapping from e-states to the events that are their images
will be used to define expected utility in the framework of evidential
structures and plans.

\Section z{Evidence-based choice}

\Subsection k{Choice alternatives, plans, and conditional preferences}

Plans, to be introduced in this section, are the entities, regarding
which conditions will be sought that ensure rationalizability in terms
of subjective conditional expected utility. A plan specifies which
alternative is chosen in each evidential state in some
evidential structure. As has already been discussed in the introduction,
the set of choice alternatives is a primitive entity in the
specification of the theory. Formally, $A$ is the set of choice
alternatives. Throughout the article, it is assumed
for mathematical simplicity that
\begin{enum} 
\meti{17}{$A$ is finite, and contains at least two distinct
  alternatives.}
\end{enum}

A function, $\zeta \from \X \to A$, is a \emph{plan.} If $\X$ is
finite, then $\zeta$ will be called a finite plan. If $\X$ is a tree,
then $\zeta$ will be called a tree plan.

A \emph{conditional-preference relation} is a ternary
relation, $\pref \subseteq A \times A \times \X$, such that
\begin{enum}
\meti{18}{$\preceq_x$ is a weak order (that is, a total preorder) on $A$, where
  $a \preceq_x b \iff (a,b,x) \in \pref$.\footnote{In all subsequent
  discussions, where reference is made to $\pref$ and $\preceq_x$,
  this relation between them is assumed. The variables, `$a\/$' and
  `$b\/$', will range over $A$.} ($\prec_x$ is the asymmetric
  part, and $\sim$ is the symmetric part, of $\preceq_x$.)}
\end{enum}

Conditional-preference relation $\pref$ \emph{rationalizes} plan
$\zeta$ iff
\display{19}{\forall x \enspace \forall a \setminus q\{ \zeta(x)
  \} \enspace a \prec_x \zeta(x)}
Plan $\zeta$ is \emph{rationalizable} if some conditional-preference relation
rationalizes it.\footnote{This terminology was introduced by
  \citet{Richter-1966}, in his abstract formulation of the economic theory of
  revealed preference. Moreover, the entire framework introduced so
  far in this article is tantamount to his framework, with an e-state
  being a \emph{budget} and a plan being a \emph{choice} in his terminology.}

\Subsection l{Immediate-dominance consistency of preferences
  and plans}

Now we introduce a putative criterion of rationality that will be a
main focus of analysis in the remainder of this paper. This principle
is the most natural cognate in the theory of evidence-based choice to
the ``sure thing principle'' in Savage's theory. Before stating the
principle, let's recall Savage's principle, beginning with the
informal example by which he motivates it.

\begin{quote}
A businessman \dots\ considers the outcome of the next presidential
election relevant to the attractiveness of a purchase. \dots\ [H]e
asks whether he would buy if he knew that the Republican candidate
would win, and decides that he would do so. Similarly, he considers
whether he would buy if he knew that the Democratic candidate
would win, and again decides that he would do so. Seeing that he would
buy in either event, he decides that he should buy, even though he
does not know which event \dots\ will
obtain[.]\footnote{\citet[chapter 2.7]{Savage-1972}.} 
\end{quote}

Keep in mind that events play the role in Savage's theory that
e-states play here and that, in fact, the non-null events of a Boolean
field, together with the subset relation, satisfy the definition of an
e-structure. Savage formulates his condition P3 as a formal postulate
corresponding to the sure thing principle, he derives its equivalence
(in the presence of other postulates) with essentially the following
assertion.\footnote{\citet[chapter 2.7, theorem
    2]{Savage-1972}. Technically, Savage's theorem is slightly
  stronger than the paraphrase of it here.} Let $\B$ be a Boolean
field of sets.

\begin{enum}
\meti{20}{\begin{quote}
If $Y \subseteq \B$ is a partition of $z \in \B$, and if a $\prec_x b$
for all $x$ in $Y$, then $a \prec_z b$.
\end{quote}}
\end{enum}

A situation in which $Y$ is the set of e-states that are immediately
more specific than $z$ (that is, $Y = \immed{z}$ according to
definition \eqn{3}), rather than a situation
in which $Y$ is a partition of $z$, will be considered here. A
conditional-preference relation is immediate-dominance (ID) consistent
iff, for every $z \in \X$,
\display{21}{ [\immed z \neq \emptyset \mand \forall x \! \in \! \immed z \; a
    \prec_x b] \implies a \prec_z b}
This condition is in a similar spirit to Savage's proposition
\eqn{20}. It expresses the idea that if, in e-state $x$, an agent
would strictly prefer $b$ to $a$ after acquiring just one piece of
further evidence beyond what is comprised in $x$, regardless of what
that new evidence would show, then the agent should already strictly
prefer $b$ to $a$ in $x$.

Call plan $\zeta$ \emph{immediate-dominance\/} (\emph{ID\/})
\emph{consistent\/} if, for every $a$ and $z$,\footnote{This
  definition corresponds to the definition of consistency in
  \citet{GreenPark-1995}, when $\zeta$ (which could be a multi-valued
  correspondence in that article) is a function.}
\display{22}{ [\immed z \neq \emptyset \mand \forall x \! \in \! \immed z \;
    [\zeta(x) = a]] \implies \zeta(z) = a}

\lema B{A plan is ID consistent if, and only if, it is rationalized
  by an ID-consistent conditional-preference relation.}

\begin{proof}
Clearly, if any ID-consistent conditional-preference relation
rationalizes plan $\zeta$, then $\zeta$ is ID consistent.  To prove
the converse, define a conditional-preference relation, $R$, by $a
\prec_x b \iff [a \neq \zeta(x) = b]$. (Otherwise, since $\preceq_x$
is a weak order, both $a \preceq_x b$ and $b \preceq_x a$.) $R$
rationalizes $\zeta$, and $\zeta$ is ID consistent if, and only if,
$R$ is ID consistent.
\end{proof}

\Section t{Representability by subjective conditional expected
  utility}

\Subsection w{Reformulating conditional probability}

The first step in examining the relationship between evidence-based
choice and maximization of conditional expected utility is to
reformulate the standard Bayesian account of conditioning beliefs on
evidence in terms of an evidential structure as defined here. Recall
the standard formulation. There is a sample space, $(\sample, \B,
\Pr)$, where $\B$ is a field of subsets of $\sample$ and $\Pr \from \B
\to [0,1]$ is a probability measure. (The field and measure may be
either finitely or countably additive.) Evidence is modeled as a
sequence, $\langle X_n \rangle_{n \in \Nat}$, of random variables,
that is, of functions $X_n \from \sample \to \Re$.\footnote{A different
  index set---for example, a finite set---may be appropriate to use in
  some specific application.}  These functions are measurable, that
is, for any $r < s$, $\{ \omega \mid r < X_n(\omega) < s \} \in \B$.

Performing potentially infinitely many tosses of a coin is an
example. A sample point (that is, an element of $\sample$), is a
function $\omega \from \Nat \to \{ 0,1 \}$, where $\omega(n) = 0$
indicates that the $n^{\text{th}}\/$ toss lands `tails' and $\omega(n)
= 1$ indicates that the $n^{\text{th}}\/$ toss lands `heads' in
possible world $\omega$. $\B$ is the Boolean field of finite
intersections of the sets (or \emph{events}), $\{ \omega \mid \omega_n =
k \}$ for $n \in \Nat$ and $k \in \{ 0,1 \}$ (or the $\sigma$ field
of countable unions and intersections of those
events). The random variable $X_n$, defined by $X_n(\omega) =
\omega(n)$, specifies how the $n^{\text{th}}\/$ toss lands in each
possible world. The epistemic states are finite sequences of
observations of toss outcomes. Again encoding `tails' by 0 and `heads'
by 1, this corresponds to
\display{23}{\X = \bigcup_{n \in \Nat} \{ 0,1 \}^n}
where $\{ 0,1 \}^0 = \{ \emptyset \}$. In possible world $\omega$, after the
coin has been tossed $n$ times, the observer is in e-state
$(\omega(0),\dotsc, \omega(n-1))$. Let $x = (x_0,\dotsc, x_{n-1}) =
(\omega(0),\dotsc, \omega(n-1))$. The event in $\B$ that is associate
with set of possible worlds in which the observer can be in this
e-state is
\display{24}{\event(x) = \bigcap_{k<n} X_k^{-1}(x_k) \in \B}

It is obvious that $\event$ satisfies analogues of conditions
\eqn{12}--\eqn{16} of \rslt C. In general, where $(\X, \wms)$ is an
e-structure and $\B$ is a field of subsets of set $\sample$, call
$\event \from \X \to \B$ an \emph{embedding} if it satisfies those
analogous conditions.\footnote{The only condition that differs from
  its analogue for $\set$, beyond a notational change, is \eqn{27}. This
  divergence is necessary because, for instance, there are no leaves
  in the e-structure that formalizes the the coin-tossing example
  discussed here. It will be proved in \sec D that, if an infinite
  e-structure is a tree, then it can be embedded in a field of
sets analogously to what has been done in \rslt C. A still more
general result is proved in \citet{bg-2016}.}  That is,
\display{25}{\event(\nothing) = \Omega}
and, for all $x$ and $y$,
\display{26}{\event(x) \neq \emptyset}
\display{27}{[\exists z \, z \sms x] \implies \event(x) =
  \bigcup_{z \sms x} \event(z)}
\display{28}{x \wms y \iff \event(x) \subseteq \event(y)}
\display{29}{x \incompat y \iff \event(x) \cap \event(y) = \emptyset}

Returning to the coin-tossing example, suppose that $\sample$ and $\B$
are as above, and that $\Pr$ is a probability measure. Then, for an
event $B \in\B$, the probability of $B$ conditional on having observed
the first $n$ tosses of the coin is
\display{30}{\Pr[B | X_0,\dotsc, X_{n-1}](\omega) = \frac{\Pr(B \cap
    \bigcap_{k<n} X_k^{-1}(x_k))}{\Pr(\bigcap_{k<n} X_k^{-1}(x_k))}}

Expressing this in the language of an e-structure and an embedding, it
becomes
\display{31}{\Pr[B | x] = \frac{\Pr(B \cap \event(x))}{\Pr(\event(x))}}

\Subsection x{Representation of conditional preferences by subjective
  conditional expected utility}

Let $\pref$ be a conditional-preference relation, for $A$ and $\X$.
Then $\pref$ is \emph{subjective-conditional-expected-utility\/}
(\emph{SCEU\/}) \emph{representable\/} if there are an embedding,
$\event$, of $(\X, \wms)$ into a measurable space, $(\sample, \B)$,
and a countably additive probability measure, $\Pr \from \B \to [0,1]$
and bounded measurable functions, $u_a \from \sample \to \Re$, for $a
\in A$, such that, for each $a$, $b$, and $x$,\footnote{All of the
  probability measures occurring in examples or in proofs of results
  in this article will be purely atomic, so $\int_{\event(x)} g \,
  d\Pr$ can be interpreted as $\sum_{x \in C} g(x) \Pr(\{ x \})$. As
  will be discussed in detail in \sec{d}, these results would continue
  to hold if countable additivity were replaced y finite additivity in
  the definition of SCEU. \label{fn:finprob}}
\display{32}{a \preceq_x b \iff \int_{\event(x)} u_a \, d\Pr \le
\int_{\event(x)} u_b \, d\Pr}
By \eqn{32}, if $\event$, $\Pr$ and $\langle u_a
\rangle_{a \in A}$ represent $\pref$, then
\display{33}{[\exists a \; \exists b \; a \prec_x b] \implies
  \Pr(\event(x)) > 0}
Note that, in view of \eqn{31}, and with conditional expectation being
defined from conditional probability by Lebesgue integration,
\display{34}{\int_B f \, d \Pr = \Pr(B) \mean[f | B]}
From \eqn{33} and \eqn{34}, it is clear that $u_a(\omega)$ should be
understood to be the state-contingent utility of alternative $a$ at
sample point $\omega$.

Regarding why boundedness of $u_a$ is required in the definition of
SCEU representability, three comments are in order.
\begin{itemize}
\item
SCEU representability is used to define \sceup ity of a plan.  \Rslt Q
will characterize \sceup\ plans on e-structures that may be infinite,
and so from which not all real-valued functions are bounded. That
result is proved by a construction of bounded functions. Including
boundedness in this definition indirectly incorporates the recognition
of that fact in the statement of the theorem.
\item
Boundedness of the utility function in Savage's representation is
implied by Savage's axioms.\footnote{\citet[theorem 14.5]{Fishburn-1970}}
Thus, incorporating boundedness in this definition would facilitate
comparisons between possible results about SCEU representability and
their analogues in Savage's framework.
\item
While some persons view decision theory as being a purely descriptive
enterprise of accounting for agents' behavior, others view it as
having a role in guiding the formulation and assessment of theories of
cognition. On this latter view, the Bayesian theory is a model of
agents whose preferences among alternativess are derived from the interaction
of two, distinct, mental entities: beliefs and state-contingent
utilities. Beliefs are represented by a probability
measure. State-contingent utilities are completely distinct from
beliefs, so they should induce a preference among acts in combination
with any possible probability measure, not only in combination with
the probability measure that the agent happens to hold. That is,
state-contingent utilities should be integrable with any probability
measure, and therefore (in order to avoid a version of the
St.\ Petersburg paradox), the state-contingent utility of any alternative
should be a bounded function.
\end{itemize}

\Subsection n{ID consistency is neither necessary nor sufficient
  for \uppercase\sceup ity}

Conditional-preference relation $\pref$ rationalizes plan $\zeta$ iff,
for all $x$, $\forall a \setminus \{ \zeta(x) \} \; a \prec_x
\zeta(x)$.

Plan $\zeta$ is \emph{\sceup\/}\ if it is
rationalized by a SCEU-representable conditional-preference relation
for $(\X,\wms)$ and $A$.\footnote{This definition corresponds closely
  to the definition of a Bayes contingent plan in
  \citet{GreenPark-1995}, except that here a plan is required to be a
  function rather than a possibly multi-valued correspondence.}
Equivalently, $\zeta$ is \sceu iff there exist an embedding,
$\event$, of $(\X, \wms)$ into a measurable space, $(\sample, \B)$, 
a probability measure, $\Pr \from
\B \to [0,1]$, and, for each $a \in A$, a bounded, measurable
function, $u_a \from \sample \to \Re$, such that, for all $x$,
\display{35}{\forall a \setminus \{ \zeta(x) \} \enspace
  \int_{\event(x)} u_a \, d\Pr < \int_{\event(x)} u_{\zeta(x)} \, d\Pr}
That is, if \eqn{35} is satisfied, then $a \prec_x b \iff
\int_{\event(x)} u_a \, d\Pr < \int_{\event(x)} u_{b} \, d\Pr$
defines an SCEU-representable conditional-preference relation that
rationalizes $\zeta$.

Note that, by \eqn{17} and \eqn{32} and \eqn{35},
\begin{enum}
\meti{36}
{if $\event$, $\Pr$ and $\langle u_a \rangle_{a \in A}$ \sceup ize
$\zeta$,\\ then $\forall x \; \Pr(\event(x)) > 0$}
\end{enum}

A more subtle question is: what is
the relation between a plan being \sceup\  and being ID
consistent? It will now be shown that, for plans in general, being ID
consistent is neither necessary nor sufficient for being SCEU
rational. It will subsequently be shown that the two conditions are
equivalent for tree plans.

Figures \ref{fig:b}--\ref{fig:d} represent plans. A label, `$x \mapsto a$', denotes
that node $x$, at that location, is mapped to alternative $a$ by the plan. In
examples \ref{xmpl:r} and \ref{xmpl:t}, the entire e-structure is the
domain of the plan.

\sbox{\labelbox}{Diagram of \Rslt r}

\begin{figure} \centering 
\begin{tikzpicture}

\draw (3,3) node[above]{$\nothing \mapsto a$} -- (2,2) node[left]{$x_1 \mapsto b$};

\draw (3,3) -- (4,2) node[right]{$x_2 \mapsto b$};

\draw (2,2) -- (1,1) node[below]{$z_1 \mapsto a$};

\draw (2,2) -- (3,1) node[below]{$z_3 \mapsto b$};

\draw (4,2) -- (5,1) node[below]{$z_2 \mapsto a$};

\draw (4,2) -- (3,1);

\end{tikzpicture}
\caption{\usebox{\labelbox}} \label{fig:b} \end{figure}

\exmpl r{ID consistency is not necessary for \sceup ity. Let $\X = \{
  \nothing, x_1, x_2, z_1, z_2, z_3 \}$. Let $\immms$ comprise
  the pairs $(x_i, \nothing)$, $(z_i, x_i)$, and $(z_3, x_i)$ (for $i
  \in \{ 1,2 \}$), and let $\wms$ be the reflexive, transitive closure of
  $\immms$. Let $A = \{ a,b \}$. Define plan $\zeta$ by
  $\zeta(\nothing) = \zeta(z_1) = \zeta(z_2) = a$ and
  $\zeta(x_1) = \zeta(x_2) = \zeta(z_3) = b$. This example is shown in
\fgr b.}

\claim s{Plan $\zeta$ of \rslt r is \sceup\ but is not ID consistent.}

\begin{proof}
Plan $\zeta$ is ID-inconsistent at $\nothing$. 

By \rslt C, $\set$ embeds $(\X, \wms)$ into $(\leaves, 2^\leaves)$,
where $\leaves =  \{ z_1, z_2, z_3 \}$. Events are
assigned by $\set$ to e-states as follows: $\set(z_i) = \{ z_i \}$,
$\set(x_i) = \{ z_i, z_3 \}$, and $\set(\nothing) = \{ z_1, z_2, z_3 \}$.

The following probability measure and utility functions constitute an
SCEU representation of a conditional-preference relation that
rationalizes $\zeta$: $\pi(z_1) = \pi(z_2) = 2/7$, $\pi(z_3) = 3/7$,
$u(a, z_1) = u(a, z_2) = u(b, z_3) = 1$, $u(b, z_1) = u(b, z_2) = u(a,
z_3) = 0$.
\end{proof}

\corol E{There exists a conditional-preference relation that is
  SCEU representable, but that rationalizes an ID-inconsistent
  conditional-preference relation.}

\begin{proof}
By claim \ref{clm:s}, there exists a SCEU-representable conditional-preference
relation that rationalizes plan $\zeta$ in \rslt r. In particular,
since the relation rationalizes $\zeta$, it satisfies $a \preceq_{x_1}
b$ and $a \preceq_{x_2} b$ but $b \preceq_{\nothing} a$. Thus it is
not ID consistent.
\end{proof}

\exmpl t{ID consistency is not sufficient for \sceup ity. Let $\X =
  \{  \nothing, x_1, x_2, z_1, z_2, z_3 \}$. Let $\immms$ comprise the
  pairs $(z_1,x_1)$, $(z_2,x_2)$, $(z_3,x_1)$, $(z_3,x_2)$,
  $(x_1,\nothing)$, and $(x_2,\nothing)$, and let
  $\wms$ be the reflexive, transitive closure of $\immms$. Let $A = \{
  a,b,c \}$. Define plan $\zeta$ by $\zeta(\nothing) = \zeta(z_3) = a$,
  $\zeta(x_1) = \zeta(z_2) = b$, and $\zeta(x_2) = \zeta(z_1) =
  c$. This example is shown in \fgr c.}

\sbox{\labelbox}{Diagram of \Rslt t}

\begin{figure} \centering 
\begin{tikzpicture}

\draw (3,3) node[above]{$\nothing \mapsto a$} -- (2,2) node[left]{$x_1 \mapsto b$};

\draw (3,3) -- (4,2) node[right]{$x_2 \mapsto c$};

\draw (2,2) -- (1,1) node[below]{$z_1 \mapsto c$};

\draw (2,2) -- (3,1) node[below]{$z_3 \mapsto a$};

\draw (4,2) -- (5,1) node[below]{$z_2 \mapsto b$};

\draw (4,2) -- (3,1);

\end{tikzpicture}
\caption{\usebox{\labelbox}} \label{fig:c} \end{figure}

\claim u{Plan $\zeta$ of \rslt t is ID consistent but is not \sceup.}

\begin{proof}
Plan $\zeta$ is ID consistent because there is no e-state, $x$, for
which $\zeta(\immed x)$ is a singleton.

Suppose that $\zeta$ were rationalized by a conditional-preference
relation represented by probability space, $(\sample, \B, \Pr)$, an
embedding, $\event \from \X \to \B$, and utility functions 
$\langle u_a \rangle_{a \in A}$.

By condition \eqn{16}, $\event(z_1) \cap \event(z_2) = \event(z_1) \cap
\event(z_3) = \event(z_2) \cap \event(z_3) =
\emptyset$. By condition \eqn{15},  $\event(x_i) = \event(z_i) \cup
\event(z_3)$ and $\event(\nothing) = \event(x_1) \cup
\event(x_2)$. It follows that
\display{37}{\sample = \event(\nothing) = \event(x_1) \cup \event(z_2) =
  \event(x_2) \cup \event(z_1)}
and
\display{38}{\event(x_1) \cap \event(z_2) =
  \event(x_2) \cap \event(z_1) = \emptyset}

Then, since $\zeta(x_1) = \zeta(z_2) = b$ and $\zeta(x_2) = \zeta(z_1)
= c$, \eqn{35} implies that $\int_{\event(x_1)} u_b \, d\Pr >
\int_{\event(x_1)} u_c \, d\Pr$ and $\int_{\event(z_2)} u_b \,
d\Pr > \int_{\event(z_2)} u_c \, d\Pr$ and $\int_{\event(x_2)}
u_c \, d\Pr > \int_{\event(x_2)} u_b \, d\Pr$ and
$\int_{\event(z_1)} u_c \, d\Pr > \int_{\event(z_1)} u_b \,
d\Pr$. These inequalities, together with \eqn{37} and \eqn{38}, imply a
contradiction: that both $\int_{\sample} u_b \, d\Pr > \int_{\sample} u_c \,
d\Pr$ and $\int_{\sample} u_c \, d\Pr > \int_{\sample} u_b \, d\Pr$.
\end{proof}

\corol F{There is a conditional-preference relation that is ID
  consistent but not SCEU representable.}

\begin{proof}
Let $R$ be the ID-consistent conditional-preference relation, guaranteed to
exist by \rslt B, that rationalizes plan $\zeta$ in \rslt t. By claim
\ref{clm:u}, $R$ is not SCEU representable (or else $\zeta$ would be \sceup).
\end{proof}

\Subsection y{A necessary condition for a finite plan to be \sceup}

In each of the proofs of claim \ref{clm:s} and claim \ref{clm:u}, \rslt C has been invoked to
translate an implication of ID consistency into the Borel-field
framework of Savage's sure-thing principle. Those translations have
revealed that ID consistency is the stronger principle in one
respect, but the weaker one in another. In \rslt r, an implication of
ID consistency translates into a condition that would be an instance
of Savage's condition \eqn{20}, except that $Y$ is only a cover of
$x$, not a partition. In that respect, ID consistency demands more
than the sure-thing principle does. But, in \rslt t, there is a
partition to which Savage's condition \eqn{20} applies, but the
partition does not correspond to a set $Y(x)$, so ID consistency has
no implication regarding it. In that respect, ID consistency demands
less than the sure-thing principle does.

ID consistency is arguably the closest analogue of the sure-thing
principle that can be formulated, without invoking a non-trivial
set-theoretic construction, for e-structures in general. However, in
the case of finite plans, a closer analogue of Savage's condition
\eqn{20} exists.\footnote{A generalization to infinite plans could be
  stated in terms of the embedding defined by equation \eqn{81}
  below.}  Its statement involves the definition of a partition of an
element of a finite e-structure.

Define $P \subseteq \X$ to be a partition of $x$ iff the following three
conditions are satisfied.
\begin{enum}
\meti{39}{If $z \in P$, then $z \wms x$}
\meti{40}{For every leaf, $z$, of $(\X, \wms)$ such that $z \wms x$,
  there is some $y \in P$ such that $z \wms y$}
\meti{41}{If $y \in P$ and $z \in P$ and $y \neq z$, then $y \incompat z$}
\end{enum}

Define a finite plan, $\zeta$, to be \emph{partition-dominance\/} (\emph{PD\/})
\emph{consistent\/} iff, for every partition, $P$, of any e-state, $x$,
\display{42}{[\forall y \! \in \! P \; \zeta(y) = a] \implies \zeta(x)
  = a}

Let $\set \from \X \to 2^\leaves$ be the embedding defined by \eqn{11}
in \rslt C. Define a finite plan, $\zeta$, to be
\emph{bilateral-dominance\/} (\emph{BD\/}) \emph{consistent\/}
consistent iff, for every pair of sets $P \subseteq \X$ and $Q
\subseteq \X$,
\display{43}{
\begin{aligned} \strut [[P \mand Q \text{\ satisfy \eqn{41}}] &\mand
{\textstyle [\bigcup_{x \in P} \set(x) = \bigcup_{x \in Q} \set(x)]}\\ \mand
[\forall x \! \in \! P \; \forall y \! \in \! P \; \zeta(x) =
    \zeta(y)] &\mand [\forall x \! \in \! Q \; \forall y \! \in \! Q \;
    \zeta(x) = \zeta(y)]]\\ \implies& [\forall x \! \in \! P \; \forall y
  \! \in \! Q \; \zeta(x) = \zeta(y)]
\end{aligned}
}

It follows immediately from \rslt C that the image under $\set$ of a
partition as defined here is a partition of an event, and that
condition \eqn{42} transforms under $\set$ into an instance of
Savage's principle \eqn{20} (translating $\zeta(x) = b$ as $a \prec_x
b$). Nonetheless, not even stronger requirement of BD consistency is
sufficient for \sceup ity, as the following example
shows.\footnote{Condition \eqn{42} follows from condition \eqn{43} by
  setting $Q = \{ x \}$. To see that \eqn{42} does not imply \eqn{43},
  add a leaf, $z_4 \immms \nothing$, to the e-structure in \rslt t and
  specify that $\zeta(z_4) = a$. No e-state of this augmented
  e-structure has a partition that satisfies the antecedent of
  \eqn{42} with respect to $\zeta$ (except, trivially, the singleton
  of the e-state itself), so $\zeta$ is PD consistent. However,
  the partitions $P = \{ x_1, z_2 \}$ and $Q = \{ x_2, z_1 \}$
  together satisfy the antecedent of \eqn{43} with respect to $\zeta$,
  but the consequent is not satisfied, so $\zeta$ is not BD
  consistent.}

\exmpl G{BD consistency is not sufficient for \sceup ity. Let
  $\X = \{ \nothing, x_1, x_2, x_3, x_4, z_1, z_2, z_3, z_4 \}$. and $A =
  \{ a_1, a_2, a_3, a_4 \}$. Define $w \wms y$ iff $w = y$ or $w =
  \nothing$ or, for $j \equiv i+1 \pmod{4}$ or $j \equiv i+2
  \pmod{4}$, $w = x_i$ and $y = z_j$. Let $\zeta(x_i) = a\zeta(z_i) =
  a_i$ and $\zeta(\nothing) \in A$ ($\zeta(\nothing)$ is irrelevant to
  the analysis of the example).  This example is shown in the left
  part of \fgr d. Note that $x \wms y$ has the geometric
  interpretation that $y$ is on a path between $\nothing$ and $x$.}

\sbox{\labelbox}{Diagram of \Rslt G}

\begin{figure} \centering 
\begin{tikzpicture}

\draw (0,4) node[above]{$z_1 \mapsto a_1$} --(2,3);

\draw (2,3.5) node{$x_4 \mapsto a_4$};

\draw (2,3) -- (4,4) node[above]{$z_2 \mapsto a_2$};

\draw (4,4) -- (3,2) node[right]{$x_1 \mapsto a_1$};

\draw (3,2) -- (4,0) node[below]{$z_3 \mapsto a_3$};

\draw (4,0) -- (2,1);

\draw (2,.5) node{$x_2 \mapsto a_2$};

\draw (2,1) -- (0,0) node[below]{$z_4 \mapsto a_4$};

\draw (0,0) -- (1,2) node[left]{$x_3 \mapsto a_3$};

\draw (1,2) -- (0,4);

\draw (2,2) node{$\nothing$};

\draw (2.64,1.42) node[rotate=-45]{$\mapsto a \in A$};

\draw (2,3) -- (2,2.3);


\draw (1,2) -- (1.7,2);

\draw (2,1) -- (2,1.7);


\draw (3,2) -- (2.3,2);

\draw (7,4) node[above]{$z_{s(1)} \mapsto a_{s(1)}$} -- (6,2) node[right]{$x_1 \mapsto a_1$};

\draw (6,2) -- (7,0) node[below]{$z_{t(1)} \mapsto a_{t(1)}$};

\end{tikzpicture}
\caption{\usebox{\labelbox}} \label{fig:d} \end{figure}

\claim S{Plan $\zeta$ of \rslt G is BD consistent but is not \sceup.}

\begin{proof}
Call $P \subseteq 2^\X$ \emph{pairwise incompatible} if it satisfies
condition \eqn{41}. If $P$ is pairwise incompatible, then call
$\bigcup_{x \in P} \set(x)$ the \emph{support} of $P$. Call $P$
\emph{monochrome} iff $\forall x \! \in \! P \; \forall y \! \in \! P
\; \zeta(x) = \zeta(y)$. In this terminology, the antecedent of
\eqn{43} states that $P$ and $Q$ are monochrome, pairwise-incompatible
sets that share the same support.

Define $s(i) \equiv i+1 \pmod{4}$ and $t(i) \equiv i+2 \pmod{4}$. The
significance of this notation is illustrated, for $i=1$, in the right
part of \fgr d.

First, it will be shown that the plan, $\zeta$, in \rslt t is BD
consistent.  Note that the only monochrome, pairwise-incompatible
subsets of $2^\X$ are $\{ \nothing \}$ and, for $1 \le i \le 4$, $\{
x_i \}$, $\{ z_i \}$, and $\{ x_i,z_i \}$. If $P$ and $Q$ share the
same support, then the cardinality of that support is $k$, where $1
\le k \le 4$. If $k=1$, then there are $i$ and $j$ such that $P = \{
z_i \}$ and $Q = \{ z_j \}$.  The two pairwise-incompatible sets share
the same support only if $z_i = \set(z_i) = \set(j) = j$, so $P=Q$. If
$k=2$, then there must be $i$ and $j$ such that $P = \{ x_i \}$ and $Q
= \{ x_j \}$.  For any $h$, the support of $\{ x_h \}$ is $\{
z_{s(h)}, z_{t(h)} \}$, and an argument parallel to the one just given
shows that $P=Q$ if the two pairwise-incompatible sets have identical
support. If $k=3$, then there are $i$ and $j$ such that $P = \{
x_i,z_i \}$ and $Q = \{ x_j,z_j \}$. For any $h$, the support of $\{
x_h,z_h \}$ is $\{ z_h,z_{s(h)},z_{t(h)} \}$, and again it is seen
that $P=Q$ if the two pairwise-incompatible sets have identical
support. If $k=4$, then $P = Q = \{ \nothing \}$. In conclusion, $P$
and $Q$ satisfy the antecedent of condition \eqn{43} only if $P=Q$. In
this trivial case, the conclusion of \eqn{43} is also satisfied, so
the condition is satisfied. That is, $\zeta$ is BD consistent.

Now it will be seen that $\zeta$ is not \sceup. A contradiction will
be derived from the assumption that probability measure $\Pr$ on
$(\sample, \B)$  and $\{ u_a \from \sample \to \Re \}_{a \in A}$
\sceup ize $\zeta$. Since the contradiction arises for all possible
choices of $\Pr$ and $\{ u_a \}_{a \in A}$, $\zeta$ is not
\sceup.

Note that, as it applies to $x_i$, condition \eqn{35} consists
of three equations, one of which specifies that $a_{s(i)} \prec_{x_i}
a_i = \zeta(x_i)$. That equation is 
\display{44}{\int_{\event(x_i)} u_{a_{s_i}} \, d \Pr <
  \int_{\event(x_i)} u_{a_i} \, d \Pr}
Defining $\q jk = \int_{\event(z_k)} u_{a_j} \, d \Pr$ and invoking
conditions \eqn{14} and \eqn{16} in the definition of an embedding,
\eqn{44} can be rewritten as
\display{45}{\q {s(i)}{s(i)} + \q {s(i)}{t(i)} <
  \q i{s(i)} + \q i{t(i)}}
Summing \eqn{45} for $i \in \{ 1,2,3,4 \}$ yields
\display{46}{\sum_{i=1}^4 \q {s(i)}{s(i)} + \sum_{i=1}^4 \q {s(i)}{t(i)} <
  \sum_{i=1}^4 \q i{s(i)} + \sum_{i=1}^4 \q i{t(i)}}
Noting that $t(i) = s(s(i))$ and $i = t(t(i))$, this can be rewritten
as
\display{47}{\sum_{i=1}^4 \q ii + \sum_{i=1}^4 \q i{s(i)} <
  \sum_{i=1}^4 \q i{s(i)} + \sum_{i=1}^4 \q {t(i)}i}
Subtracting $\sum_{i=1}^4 \q i{s(i)}$ from both sides yields
\display{48}{\sum_{i=1}^4 \q ii < \sum_{i=1}^4 \q {t(i)}i}
However, for each $i$, since $\zeta(z_i) = a_i$, $\q {t(i)}i < \q
ii$. Summing this inequality over $i$ yields
\display{49}{\sum_{i=1}^4 \q {t(i)}i < \sum_{i=1}^4 \q ii}
which contradicts \eqn{48}.
\end{proof}

\Subsection A{ID consistency is equivalent to \sceup ity for finite tree plans}

An e-structure, $(\X, \wms)$, is a \emph{tree} if the following
condition is satisfied for every e-state, $x$.\footnote{This
  definition and the following implication are easily shown to be
  equivalent to common definitions of an arborescence (that is, a
  directed tree), such as \citet[definition 3.20]{Gallier-2011},
  except that, unlike some of those definitions, the present
  definition does not require $\X$ to be finite.}
\display{50}{ \{ y \mid x \wms y \} \text{\ is finite, and\ }
\forall y \, \forall z \, [[x \wms y \mand x \wms z] \implies
[y \wms z \text{\ or\ } z \wms y]]}

A \emph{path\/} of length $n$ from $x$ to $y$ is specified by a
sequence $\langle z_0,\dotsc, z_n \rangle$, that satisfies the
following condition.
\display{51}{\begin{split}
z_0 = x& \mand z_n = y \\
\forall k \! < \!n& \; z_k \immms z_{k+1} \qquad \text{if\ } x \sms y\\
\forall k \! < \!n& \; z_{k+1} \immms z_k \qquad \text{if\ } y \sms x
\end{split}}
If $y$ is the root, then $n$ is defined to be is the \emph{depth\/} of
$x$.\footnote{It follows from the uniqueness assertion in \rslt J that
  there is a unique number satisfying the definition of the depth of
  an e-state.}

The following lemma summarizes standard facts regarding
a tree. It follows immediately from \citet[theorem 3.2(6)
  and theorem 3.3(3)]{Gallier-2011}.

\lema J{Suppose that $(\X,\wms)$ is a tree.  For all $x$ and $y$, $x
  \sms y$ if, and only if, there is a path from $x$ to $y$ satisfying
  $\forall k \! < \!n \; z_k \immms z_{k+1}$. There is a path from $x$
  to $y$ if, and only if, there is a path from $y$ to $x$. If there is
  a path from $x$ to $y$, then it is unique.  If $x \sms y$, then the
  depth of $x$ is strictly greater than the depth of $y$.}

The following proposition is the main result of this section.

\propo H{A finite tree plan is \sceu if, and only if, it is ID consistent.}

The proof proceeds via several lemmas. \Rslt P shows that, if an
e-structure is a tree, then an embedding associates $Y(x)$ with a
partition of the image of $x$ under the embedding.  \Rslt I concerns
induction and recursion on well founded relations. Those principles
are used to prove \rslt K, showing that for any \sceu plan, $\zeta$,
e-state, $x$, and alternative, $a \neq \zeta(x)$, it is possible to
``drill down'' to a leaf that is weakly more specific than $x$,
constructing the path from $x$ to the leaf in such a way that $\zeta$
never specifies $a$. Lemmas \ref{lem:M} and \ref{lem:N} establish further,
technical details that will be used in the main proof.

\lema P{If an e-structure, $(\X, \wms)$, is a tree; and if $\event$
  embeds $(\X, \wms)$ into a Boolean field of sets, $(\Omega, \B)$;
  and if $\exists z \, z \sms x$; then $\{ \event(y) \mid
  y \in Y(x) \}$ is a partition of $\event(x)$.}

\begin{proof}
It must be shown that $\{ \event(y) \mid y \in Y(x) \}$ is a union of
nonempty, pairwise disjoint sets, the union over which is $\event(x)$.
By the definition of an embedding, for every $y$, $\event(y)
\neq \emptyset$. If $w \neq y$ and $\{ w,y \} \subseteq Y(x)$, then
not $w \wms y$ and not $y \wms w$, by the antisymmetry of $\wms$ and
the definition of $\immms$. Then, by the definition of a tree, there
is no $z$ such that $w \wms z$ and $y \wms z$. That is, $w \incompat
y$. Therefore, by the definition of an embedding, $\event(w) \cap
\event(y) = \emptyset$. It remains only to show that $\event(x) =
\bigcup_{y \in Y(x)} \event(y)$. By the definition of an embedding,
$\event(x) = \bigcup_{z \sms x} \event(z)$. Suppose that $z \sms
x$. Then, either $z \in Y(x)$ or else, for some $y$, $z \sms y \sms
x$. Since $(\X, \wms)$ is tree, $\{ y \mid z \sms y \sms x \}$ is finite,
so there is actually some $y$ such that $z \sms y \immms x$. By the
definition of an embedding, $\event(z) \subseteq \event(y)$. Thus,
$\event(x) = \bigcup_{z \sms x} \event(y)$ implies that $\event(x) =
\bigcup_{y \in Y(x)} \event(y)$.
\end{proof}

In the following lemma, proved in \citet[theorem 5.8]{Gallier-2011},
the notation, `$f \restrict C\/$', denotes the restriction of function
$f$ to subset $C$ of its domain.

\lema I{Let $\prec$ be a well founded relation on a set, $X$. If $W \subseteq
  X$ is such that $\forall x \; [[\forall w \! \prec \! x \; w \in W]
    \implies x \in W]$, then $W=X$. If, where $Z$ is an arbitrary,
  non-empty set, $F = \{ f \mid f \from W \to Z \mand W
  \subseteq X \}$~and $g \from F \times X \to X$, then there is a
  unique function, $f \from X \to Z$, such that $\forall x \; f(x) =
  g(f \restrict \{ z \mid z \prec x \}, x)$.}

\lema K{Suppose that $(\X, \wms)$ is a finite e-structure that is a
  tree, and that $\zeta \from \X \to A$ is ID consistent. Since $\X$
  is finite, it may be assumed without loss of generality that $\X =
  \{0,1,\dotsc,n \}$ for some number, $n$. Define the set of ``bad
  choices'' by
\display{52}{\badchoice = \{ (x,a) \mid a \neq \zeta(x) \}}
There is a unique function, $\avoid \from \badchoice \to
\leaves$, that satisfies the following definition.\footnote{Recall
  definition \eqn{10}, that $\leaves$ is the set of leaves of $\wms$.}
\display{53}{\avoid(x,a) = \begin{cases}
x &\text{if\ } x \in \leaves\\
\min \{ \avoid(y,a) \mid (y,a) \in (Y(x) \times \{ a \}) \cap \badchoice \},
&\text{if\ } x \notin \leaves
\end{cases}}
Function $\avoid$ is surjective and, for all $(x,a)$ in $\badchoice$,
satisfies the condition that
\display{54}{\avoid(x,a) \wms x \mand \forall y \; [[\avoid(x,a) \wms y
      \wms x] \implies \zeta(y) \neq a]}
}

\begin{proof}
Suppose that $x \notin \leaves$ and $(x,a) \in \badchoice$. First,
note that $(Y(x) \times \{ a \}) \cap \badchoice \neq
\emptyset$. Otherwise, ID consistency would require that $\zeta(x) =
a$, contradicting $(x,a) \in \badchoice$. Thus the minimum in the
clause of \eqn{53} that defines $\avoid(x,a)$ is a well defined
natural number. It needs to be shown that this number is an element of
$\leaves$. This follows by well founded induction on $\wms$, which is
well founded because $\X$ is finite. To apply \rslt I, take $W = \{ y
\mid \forall a \; [(y,a) \in \badchoice \implies \avoid(y,a) \in
  \leaves] \}$.

$\avoid$ is surjective because, for every leaf, $x$, there is an
alternative, $a$, such that $\zeta(x) \neq a$; and $\avoid(x,a) = x$.

Assertion \eqn{54} also holds by well founded induction on $\wms$. This will be
proved separately for the two conjuncts of the assertion.  To see that
$\forall x \; \forall a \; [(x,a) \in \badchoice \implies \avoid(x,a)
  \wms x]$, take $X = \X$ in \rslt I, and define
\display{55}{\begin{split}
W = \{ w\mid \forall a \; [[(w,a) \in \badchoice \implies
\avoid(w,a) \wms w] \mand\\ \forall y \; [[\avoid(w,a) \wms y
      \wms w] \implies \zeta(y) \neq a]] \} 
\end{split}}
\Rslt I asserts that if, for all $x$, $x \in W$ is implied by $\forall
y \! \sms \! x \; y \in W$. In fact, it will be proved here that, for
all $x$, $x \in W$ is implied by $\forall \, y \immms x \;\strut\, y
\in W$, so \eqn{54} will follow from \rslt I. 
Consider first the induction hypothesis that, for all $w \in
Y(X)$ and for all $a$, $(w,a) \in \badchoice \implies \avoid(w,a) \wms
w$, and suppose that $(x,b) \in \badchoice$. If $x$ is a leaf and
$(x,b) \in \badchoice$, then condition \eqn{53} specifies that
$\avoid(x,b) = x \wms x$. Since $b$ is arbitrary, this shows that
$\forall a \; [[(x,a) \in \badchoice \implies \avoid(x,a) \wms w]$. If
  $x$ is not a leaf, then let $y = \min \{ \avoid(w,b) \mid (w,b) \in
  (Y(x) \times \{ b \}) \cap \badchoice \}$. Condition \eqn{53},
  together with the induction hypothesis, specifies that $\avoid(x,b)
  = \avoid(w,b) \wms w \wms x$. Again, since $b$ is arbitrary,
  $\forall a \; [[(x,a) \in \badchoice \implies \avoid(x,a) \wms
      x]$. Since this condition holds whether or not $x$ is a leaf, it
    holds for all $x$.

To prove the second conjunct of \eqn{54} by well founded induction on
$\wms$, consider the induction hypothesis that, for all $w \in Y(X)$
and for all $a$ and $y$, $[(w,a) \in \badchoice \mand \avoid(w,a) \wms
  y \wms w]\\ \implies \zeta(y) \neq a$. Suppose that $(x,b) \in
\badchoice$ and, for some $z$, $\avoid(x,b) \wms z \wms x$. If $x$ is
a leaf, then $z = x$, so $(x,b) \in \badchoice$ implies that $\zeta(z)
\neq b$.  If $x$ is not a leaf, then let $y = \min \{ \avoid(w,b) \mid
(w,b) \in (Y(x) \times \{ b \}) \cap \badchoice \}$, so that
$\avoid(x,b) = \avoid(y,b)$, so the first assertion of \eqn{54}
(proved in the preceding paragraph) entails that $\avoid(x,b) \wms y$
Thus, by \rslt J, there is a path from $\avoid(x,b)$ to $y$, which
extends to $x$ since $y \immms x$. Also by \rslt J (applying it twice
and concatenating paths, if $\avoid(x,b)$ and $z$ and $x$ are all
distinct), since $\avoid(x,b) \wms z \wms x$, there is a path from
$\avoid(x,b)$ to $x$ that passes through $z$. Since the path from
$\avoid(x,b)$ to $x$ is unique, the two paths just enumerated are
identical. Thus, either $\avoid(x,z) = \avoid(y,b) \wms z \wms y$ or
else $z = x$. In the former case, the induction hypothesis entails
that $\zeta(z) \neq b$. In the latter case, since $(x,b) \in
\badchoice$, $\zeta(z) = \zeta(x) \neq b$. The overall conclusion from
these various cases is that, for all $a$ and $y$, $[(x,a) \in
  \badchoice \mand \avoid(x,a) \wms y \wms x] \implies \zeta(y) \neq
a$, so assertion \eqn{54} follows from \rslt I (that is, by well
founded induction).
\end{proof}

Recall the definition \eqn{10}, that $\set(x) = \{ z \mid z \in \leaves
\mand z \wms x \}$. The following characterization of $\set$ follows
easily from \rslt K.

\corol O{If the e-structure, $(\X, \wms)$, is a finite tree, then\\ $z
  \in \set(x) \iff \exists y \; [y \wms x \mand \exists a \,
    z = \avoid(y,a)]$.}

The next lemma will be used to construct a \sceu state-contingent
utility functions, in terms of which an ID-consistent plan can be
rationalized. Given a set, $\range$, and an embedding, $\br$, of $(\X,
\wms)$ in $(\range, 2^\range)$ and a surjective function, $\atom$,
from some set, $\domain$, to $\range$, a new embedding, $\event$, is
defined that maps each e-state to a subset of the graph, $\graph$, of
$\atom$. That is, $\event \from \X \to 2^\graph$ and
\display{56}{\graph = \{ (b, s) \mid b \in \domain \mand s
  \in \range \mand \atom(b) = s \}}
The lemma is proved by verifying routinely that each of the conditions
\eqn{25}--\eqn{29} in the definition of an embedding is inherited by
$\event$ from $\br$. The lemma applies with $\domain = \badchoice$,
$\range = \leaves$, $\br = \set$, and $\atom = \avoid$ by \rslt C and the
surjectivity assertion of \rslt K.

\lema M{If $\br$ embeds $(\X, \wms)$ in $(\range, 2^\range)$, $\atom$
  maps an arbitrary set, $\domain$, surjectively to $\range$,
  $\graph$ is defined by \eqn{56}, and $\event$ is defined as follows by
  \eqn{57}, 
\display{57}{\event(x) = \{ (b, s) \mid \atom(b) = s \in
  \br(x) \} = \graph \cap (\domain \times \br(x))}
then $\event$ embeds $(\X, \wms)$ in $(\graph, 2^\graph)$.}

The final prerequisite to proving \rslt H is to assign natural
numbers to the elements of $\X \times A$ in a way that respects $\sms$
as an order relation.

\lema N{If $(\X,\wms)$ is a tree, and if $\X$ is countable (but not
  necessarily finite) and $A$ is finite, then there is an injective
  function $\num \from \X \times A \to \Nat$ such that, for $(x,a) \in
  \X \times A$ and $(y,b) \in \X \times A$,
\display{58}{x \sms y \implies \num(x,a) > \num(y,b)}
}

\begin{proof}
Assume, without loss of generality, that $\X \subseteq \Nat$ and $A =
\{ 0,\dotsc,m \}$.  Let $k_0,k_1\dots$ enumerate without repetition
the prime numbers greater than $2^m$.\footnote{There are infinitely
  many such primes, so distinct elements of $\X$ can be associated
  with distinct primes. Cf.~\citet[theorem 5.6]{Gallier-2011}. Also
  invoked below is the unique prime factorization theorem, that each
  $n > 1$ has a unique (up to permutation) factorization into powers
  of prime numbers. Cf.~\citet[theorem 5.10]{Gallier-2011}.} For
$(x,a) \in \badchoice$, define $\num(x,a) = 2^a \prod_{\{ z \mid x
  \wms z \}} k_z$. (Recall that $a$ and $z$, which are respectively an
alternative and an e-state, are being represented as natural numbers.)

To prove that $\num$ is injective, note that if $(x,a) \neq (y,b)$,
then either $a \neq b$ or else $x \neq y$, which is equivalent to not
both $x \wms y$ and $y \wms x$ because $\wms$ is reflexive and
antisymmetric. It must be shown that, in either case, $\num(x,a) \neq
\num(y,b)$. Consider first the case that $a \neq b$. Then there are
odd numbers, $i$ and $j$, such that $\num(x,a) = 2^a i$ and $\num(y,b)
= 2^b j$.  By the unique-prime-factorization theorem, $2^a i \neq 2^b
j$.  That is, $\num(x,a) \neq \num(y,b)$. Next consider the case that
$x \neq y$. Suppose, without loss of generality, that not $y \wms
x$. Then, $k_y$ is not a prime factor of $\num(x)$, but it is a prime
factor of $\num(y)$. Again by the unique-prime-factorization theorem
$\num(x,a) \neq \num(y,b)$. Since $\num(x,a) \neq \num(y,b)$ is
obtained in both cases, $\num$ is injective.

Suppose now that $x \sms y$. By \rslt J, $\depth(y) < \depth(x)$, and
also $\num(y,b) = 2^b \prod_{\{ z \mid y \wms z \}} k_z$ and $\num(x,a) =
2^a \prod_{\{ z \mid x \wms z \}} k_z = 2^{(a-b)} \prod_{\{ z \mid x \wms
  z \sms y \}} k_z \num(y,b) \ge 2^{\mathord{-}b} k_x \num(y,b)$.
Since $b\le m$, $2^b < k_x$, so $\num(y,b) < \num(x,a)$.
\end{proof}

Now the groundwork has been laid to prove \rslt H, stating that a
finite tree plan is \sceup\ if, and only if, it is ID consistent.  For
the remainder of this section, the conditions of \rslt K are assumed
to hold. That is, $(\X, \wms)$ is a finite e-structure that is a tree,
and $\zeta \from \X \to A$.  $\badchoice$ and
$\avoid$ are defined by equations \eqn{52} and \eqn{53} respectively.
The convention in the proof of \rslt N, that $\X$ and $A$ are
represented as subsets of $\Nat$, continues to be observed.

\begin{proof}[Proof of \Rslt H]
First it is proved that, if $\zeta$ is ID consistent, then it is
\sceup. A sample space, $(\sample[\zeta], 2^{\sample[\zeta]}, \Pr)$,
will be employed, where $\sample[\zeta]$ is essentially the graph of
$\avoid$.\footnote{An element of $\sample[\zeta]$ is of form
  $(x,a,\avoid(x,a))$, while an element of the graph of $\avoid$ is of
  form $((x,a),\avoid(x,a))$.}
\display{59}{\sample[\zeta] = \{ (x,a,\avoid(x,a)) \mid (x,a) \in \badchoice \} }
The elements of $\sample[\zeta]$ will be called \emph{sample
  points}.\footnote{Sample points are what Bayesian
  theorists, following \citet{Savage-1972}, generally call
  \emph{states of the world}.}

The probability measure, $\Pr$, is now defined by specifying
its values on the singleton events in $2^{\sample[\zeta]}$. Since $\sample[\zeta]$ is
finite, the probabilities of all larger events are determined by additivity.
\display{60}{\begin{aligned}
q_{(x,a,z)} &= \num(x,a)\\
r_\omega &= 3^{\mathord{-}q_\omega}\\
\Pr( \{ \omega \} ) &= \frac{r_\omega}{\sum_{\psi \in \sample[\zeta]} r_\psi}
\end{aligned}}
Note that, since there is a 1--1 correspondence between the domain of
a function ($\avoid$, in this case) and its graph, and since $\num$ is
an injective function, $\psi \neq \omega$ implies that $q_\psi \neq
q_\omega$. 

State-contingent utility functions, $\{ u_a \from \sample[\zeta] \to \Re \}_{a
  \in A}$ will now be defined.
\display{61}{u_a(w,b,z) = \begin{cases}
1, & \text{if\ } \exists y \; [z \wms y \wms w \mand \zeta(y) = a]\\
0, & \text{otherwise}
\end{cases}}

Define $\event[\zeta] \from \X \to 2^{\sample[\zeta]}$ by
\display{62}{\event[\zeta](x) = (\X \times A \times \set(x)) \cap \sample[\zeta]}
By \rslt M, $\event[\zeta]$ is an embedding.

According to condition \eqn{35}, with integration expressed as summation
as specified in footnote \ref{fn:finprob}, $\event[\zeta]$, $\Pr$, and $\{
u_a \from \sample[\zeta] \to \Re \}_{a \in A}$ \sceup ize $\zeta$ iff, for
every $(x,a)$ in $\badchoice$,
\display{63}{\sum_{\omega \in \event[\zeta](x)} [u_{\zeta(x)}(\omega) -
    u_a(\omega)] \Pr( \{ \omega \} ) > 0}
If $(y,a,z) \in \event[\zeta](x) = (\X \times A \times \set(x)) \cap
\sample[\zeta]$, then $z \wms x$ by \eqn{11}, and also $z \wms y$ by \eqn{54}
and \eqn{59}. Then, by the definition of a tree, \eqn{50}, 
\display{64}{(y,a,z) \in \event[\zeta](x) \implies [x \wms y \text{\ or\ } y \wms x]}

The next step is to fix $(x,a) \in \badchoice$ and to express \eqn{63}
in terms of a partition of $\event[\zeta](x)$. Let $\omega =
(x,a,\avoid(x,a))$. By \eqn{59} and \eqn{62}, $\omega \in
\event[\zeta](x)$. Denote the set of sample points, $(y,b,z)$, in $\event[\zeta](x)$
that are as shallow as $\omega$ (in the sense that $x \wms y$),
excluding $\omega$, by $S$, and the set of sample points in $\event[\zeta](x)$ that
are deeper than $x$ by $D$. These sets are defined formally by
\display{65}{S = \sample[\zeta] \cap \{ (y,b,z) \mid z \wms x \wms y \}
  \setminus \{ \omega \} \qquad D = \sample[\zeta] \cap \{ (y,b,z) \mid z
  \wms y \sms x \}}
By \cor O, $\{ \{ \omega \}, S, D \}$ is a partition of
$\event[\zeta](x)$. Therefore \eqn{63} can be written as
\display{66}{\begin{split}
[u_{\zeta(x)}(\omega) - u_a(\omega)] \Pr(\omega)
+ \sum_{\psi \in S} [u_{\zeta(x)}(\psi) - u_a(\psi)] \Pr( \{ \psi \} )\\
+ \sum_{\psi \in D} [u_{\zeta(x)}(\psi) - u_a(\psi)] \Pr( \{ \psi \} ) > 0
\end{split}}

The terms of of \eqn{66} corresponding to $\omega$, $S$, and $D$ are
now examined, beginning with\\ $[u_{\zeta(x)}(\omega) - u_a(\omega)]
\Pr(\omega)$.

By \eqn{54}, $\avoid(x,a) \wms x$. Therefore $u_{\zeta(x)}(\avoid(x,a))
= 1$ by \eqn{61}, the definition of $u_a$. Also by \eqn{54} and that
definition, $u_a(\avoid(x,a)) = 0$. Therefore
\display{67}{[u_{\zeta(x)}(\omega) - u_a(\omega)]
\Pr(\omega) =  \Pr(\omega)}

Next, consider $\sum_{\psi \in S} [u_{\zeta(x)}(\psi) - u_a(\psi)]
\Pr( \{ \psi \} )$ and $\sum_{\psi \in D} [u_{\zeta(x)}(\psi) -
  u_a(\psi)] \Pr( \{ \psi \} )$. Suppose that $\psi =
(y,b,z)$.  Condition \eqn{65}, defining
$\psi \in S$, is that $z \wms x \wms y$ and $\psi \neq \omega$, so
$u_{\zeta(x)}(\psi) = 1$ by \eqn{61}. Since $u_a(\psi) \le 1$,
\display{68}{\sum_{\psi \in S} [u_{\zeta(x)}(\psi) - u_a(\psi)] \Pr( \{ \psi \} ) > 0}
Also, because $u_a(\psi) \le 1$ and $u_{\zeta(x)}(\psi) \ge 0$
\display{69}{\sum_{\psi \in D} [u_{\zeta(x)}(\psi) - u_a(\psi)] \Pr( \{ \psi \} ) \ge
  \mathord{-}\sum_{\psi \in D} \Pr( \{ \psi \} )}

It follows from \eqn{67}--\eqn{69} that 
\display{70}{\begin{split}
[u_{\zeta(x)}(\omega) - u_a(\omega)] \Pr(\omega)
+ \sum_{\psi \in S} [u_{\zeta(x)}(\psi) - u_a(\psi)] \Pr( \{ \psi \} )\\
+ \sum_{\psi \in D} [u_{\zeta(x)}(\psi) - u_a(\psi)] \Pr( \{ \psi \} ) > \Pr( \{
\omega \} - \sum_{\psi \in D} \Pr( \{ \psi \} )
\end{split}}
Thus, condition \eqn{63} is a sufficient condition for ID consistency
to imply \sceup ity, and a sufficient condition for \eqn{63} is that
\display{71}{\sum_{\psi \in D} \Pr( \{ \psi \} ) < \Pr( \{ \omega \})}
By condition \eqn{60}, $q_\omega = \num(x,a)$. Also by conditions
\eqn{58}, \eqn{60}, and \eqn{65}, $q_\psi > \num(x,a)$ for $\psi \in
D$. Therefore, by \eqn{60}, $r_\omega \ge 3^{\mathord{-}\num(x,a)}$ and,
for $\psi \in D$, $r_\psi < 3^{\mathord{-}\num(x,a)}$. Moreover, as was
remarked just after definition \eqn{60} of $q$, there is a 1--1
correspondence between $\psi$ and $q_\psi$. Thus
\display{72}{
\sum_{\psi \in D} r_\psi = \sum_{\psi \in D} 3^{\mathord{-}q_\psi}
\le \sum_{n > \num((x,a)} 3^{\mathord{-}n} < 3^{\mathord{-}\num(x,a)}
= 3^{\mathord{-}q_{\omega}} = r_{\omega} }
Since $\Pr( \{ \xi \} )$ is a positive-real multiple of $r_\xi$ (by a
factor that does not depend on $\xi$), \eqn{72} implies \eqn{71}, so
being ID consistent is sufficient for a plan to be \sceup.

Finally the converse implication, from \sceup ity of plan $\zeta$ to ID
consistency, must be proved. Suppose that $\event \from \X \to \sample$
is an embedding into a probability space, $(\sample, \B, \Pr)$ that,
along with state-contingent utility functions, $\langle u_a \rangle_{a
  \in A}$, \sceup izes $\zeta$. Suppose that $x \notin \leaves$ and
that, for all $y \in Y(x)$, $\zeta(y) = a$. Let $b \neq a$. Then,
since $\zeta$ is \sceup, for all $y \in Y(x)$, $\int_{\event(y)} u_a
- u_b \, d \Pr > 0$. By \rslt P, $\{ \event(y) \mid y \in Y(x) \}$
is a partition of $\event(x)$. Therefore, $\int_{\event(x)} u_a
- u_b \, d \Pr > 0$. From this inequality holding for all $b \neq a$,
it follows from \sceup ity that $\zeta$ is ID consistent at $x$.
\end{proof}

\Section B{Countably branching tree plans}

\Subsection D{Why study infinite plans?}

There are some questions that arise specifically in the context of
infinite e-structures and plans. For example, question, under what
conditions do an agent's conditional beliefs about an event converge
asymptotically to certainty, can only be studied if the agent might
proceed through an infinite sequence of successively more specific
evidential states. That study requires an e-structure resembling the
one contemplated in \sec w with respect to betting on the bias of a
coin that is tossed potentially infinitely many times, that consisted
of of all finite binary sequences. Besides there being applications
that require an infinitely deep evidential structure to model, there
are also some applications that require an infinitely broad one. If,
in some e-state, $x$, a chemist receives a reading of how many clicks
of a Geiger-Mueller counter were registered during a one-hour interval
in the presence of some sample of slightly radioactive material. The
number of clicks might be any non-negative integer, so it might be
argued that, in principle, $Y(x)$ should be a countable
set.\footnote{Of course, an argument about this point has two
  sides. Offsetting the desirability of adopting a specification that
  imposes no ad hoc bound, is the consideration that any physical
  measuring apparatus will have a finite upper bound on how many
  distinct events can be counted during an hour.}  It turns out that
countably branching e-structures are no more difficult to analyze than
finitely-branching ones are, so, in this section, it will be assumed
that
\display{73}{\text{For every\ } x, \, Y(x) \text{\ is a countable set.}} 
Obviously there are cases (for instance, in which evidence consists of
real-valued measurements) in which, in principle, $Y(x)$ ought to be
uncountable. That extension of the theory will not be made here,
though, because it would entail a big sacrifice of mathematical
tractability.

\Subsection E{Equivalence of \sceup ity and ID consistency}

The following generalization of \rslt H is to be proved.

\propo Q{If $(\X, \wms)$ is both an e-structure and a countably
  branching tree, and if $\zeta \from \X \to A$ is ID consistent, then
  $\zeta$ is \sceup. Conversely, under those assumptions, a \sceu tree
  plan must be ID consistent.}

The proof, as well as the statement, of this proposition follows
closely its finite analogue. Before pursuing the details of this
proof, let us take stock of what has to be done in order to dispense
with finiteness. The proof of \rslt H depends on \rslt C, lemmas
\ref{lem:J}--\ref{lem:N}, and \rslt O. Among these, lemmas
\ref{lem:J}--\ref{lem:I} make no mention of cardinality, and so can be
used in the proof of \rslt Q. \Rslt N requires that $\X$ should be
countable, but not necessarily finite. It will be proved, as \rslt R,
that a countably branching tree has countably many nodes, so \rslt N
can also be used here.

Proofs of the remaining results on which \rslt H depends, and the
proof of the proposition itself, depend on finiteness in two
ways. First, in the construction of a probability space,
$(\sample[\zeta], \B, \Pr)$, to represent $\zeta$ as a
\sceup\ rational plan, $\sample[\zeta]$ is taken to be a set that is
in 1--1 correspondence with $\leaves$, and the fact that $\set$ embeds
$(\X, \wms)$ in $(\leaves, 2^\leaves)$ (that is, \rslt C) is used to
construct an embedding, $\event$, of $(\X, \wms)$ in $(\sample[\zeta],
\B)$. However, as has just been pointed out, an infinite tree may have
no leaves, so the usefulness of \rslt C is restricted to the finite
case. To prove \rslt Q, a generally valid embedding must be devised to
play the role that \rslt{M} plays in the finite case. This will be
done in \rslt{L}. Second, \rslt I is proved by well founded recursion
and induction on $\wms$, and there are some infinite trees for which
$\wms$ is not well founded.\footnote{For example, the e-structure in
  the coin-tossing example of \sec{w} is not well founded.}
Specifically, $\avoid$ is defined recursively on $\wms$ in the proof
of \rslt{I}, and a substitute for that construction must be
found. That will be done in \rslt{A}, where recursion on the natural
numbers will be used to construct a different function (also to be
called $\avoid$) with domain $\badchoice$, and then to have the range
of $\avoid$ play the role that $\leaves$ plays in the finite case.

\lema R{If $(\X, \wms)$ is a countably branching tree, then $\X$ is countable.}

\begin{proof}
Let $(\X, \wms)$ be a countably branching tree. By the definition of a
tree, each element of $\X$ has a finite depth. Define $\X_n = \{ x
\mid x \in \X \text{\ and the depth}\\ \text{of\ } x \text{\ is\ } n \}$. By
induction, each $\X_n$ is is countable. The basis case is that $\X_0 =
\{ \nothing \}$. Suppose that $\X_n$ is countable. Then $\X_{n+1} =
\bigcup_{z \in \X_n} Y(z)$ The right side of this equation is a
countable union of countable sets, which is countable, so $X_{n+1}$ is
countable. Therefore, by induction, every $\X_n$ is
countable. Finally, $\X = \bigcup_{n \in \Nat} \X_n$ is a countable
union of countable sets, so it is countable.
\end{proof}

\lema L{Suppose that $(\X, \wms)$ is an e-structure and $\range$ is an
  arbitrary set, and that $\domain \subseteq \X \times A$ and $\atom
  \from \domain \to \range$ satisfy the following conditions.
\display{74}{\atom \from \domain \to \range
\text{\ is surjective, and\ } \forall x \; \exists a \; (x,a) \in
\domain}
\display{75}{\begin{aligned}
\forall x \; \forall a \; [[(x,a) \in \domain
      \mand \exists y \, y \sms x]& \implies\\ \exists z \; [z \sms x 
      \mand (z, a) \in \domain &\mand \atom(z,a) =
      \atom(x,a)]]
\end{aligned}}
\display{76}{\forall x \; \forall y \; [x \incompat y \implies
    \atom(\domain \cap (\{ x \} \times A)) \cap \atom(\domain \cap (\{
    y \} \times A)) = \emptyset]}
Define $\br \from \X \to 2^{\range}$ by
\display{77}{\br(z) = \atom(\domain \cap (\{ w \mid w \wms z \} \times A))}
Then $\br$ is an embedding of $(\X, \wms)$ into $(\range, 2^{\range})$.}

\begin{proof}
It must be proved that $\br$ satisfies conditions \eqn{25}--\eqn{29}
of the definition of an embedding. By \eqn{77}, $\br(\nothing) = 
\atom(\domain \cap (\{ w \mid w \wms \nothing \} \times A)) =
\atom(\domain)$, so $\br(\nothing) = \range$ because, by \eqn{74},
$\atom$ is surjective. That is, \eqn{25} holds.

Condition \eqn{26}, that $\br(x) \neq \emptyset$, holds by definition
\eqn{77} and because, by \eqn{74}, $\domain \cap (\{ x \} \times A)$
is not empty.

By \eqn{77}, $\bigcup_{z \sms x} \br(z) \subseteq \br(x)$. To verify
\eqn{27}, it remains to show that $\br(x) \subseteq \bigcup_{z \sms x}
\br(z)$ if $\exists y \, y \sms x$. Suppose that $\exists y \, y \sms
x$ and $s \in \br(x)$. By \eqn{75}, there are some $y \wms x$ and $a
\in A$ such that $(y,a) \in \domain$ and $\atom(y,a) = s$. If $y \sms
x$, then $s \in \bigcup_{w \sms x} \br(w)$. Otherwise, $y = x$. By
\eqn{75}, there is some $z \sms x$ such that $(z,a) \in \domain$ and
$\atom(z,a) = s$. Again, $s \in \bigcup_{w \sms x} \br(w)$.

Condition \eqn{28} follows directly from definition \eqn{77} of $\br$.

To prove \eqn{29}, first suppose that $x \incompat y$. Let $w \wms x$
and $z \wms y$. Then $w \incompat z$. By \eqn{76}, $\atom(\domain \cap
(\{ w \} \times A)) \cap \atom(\domain \cap (\{ z \} \times A)) =
\emptyset]$. Therefore, since this conclusion holds for all $w$ and
  $z$ that are weakly more specific than $x$ and $y$ respectively,
  $\br(x) \cap \br(y) = \emptyset$. Conversely, suppose that $\br(x)
  \cap \br(y) = \emptyset$. Then there do not exist $z$ and $a$ such
  that $z \wms x$ and $z \wms y$ and $(z,a) \in \domain$, or else
  $\atom(z,a) \in \br(x) \cap \br(y)$. But if not $x \incompat y$,
  then there is some $z$ such that $z \wms x$ and $z \wms y$, in which
  case also \eqn{74} would guarantee that $(z,a) \in \domain$ for some
  $a$. Thus it is proved by contradiction that $x \incompat y$,
  completing the proof of \eqn{29}.
\end{proof}

Now an analogue of \rslt K for tree plans, which need not be finite, is
stated and proved.

\lema A{Suppose that $(\X, \wms)$ is an e-structure that is a tree,
  and that $\zeta \from \X \to A$ is ID consistent. By \rslt R, it may
  be assumed without loss of generality that $\X \subseteq \Nat$.
  Define $\badchoice$ by \eqn{52}.  Then there is a unique function,
  $\avoid \from \badchoice \to \X^\Nat$ such that, for all $(x,a) \in
  \badchoice$, $\avoid(x,a)$ is the function $f$, that is defined
  according to \eqn{78}. Let $\langle z_0, \dotsc, z_n \rangle$ be the
  path from $\nothing$ to $x$ (that is, satisfy \eqn{51} with $z_0 =
  \nothing$ and $z_n = x$).
\display{78}{f(k) = \begin{cases}
z_k & \text{if\ } k \le n\\
\min \{ y \mid y \in Y(f(k-1)) \mand \zeta(y) \neq a \} &\text{if\ } k > n \mand\\
& \hspace{-1ex} Y(f(k-1)) \neq \emptyset\\
f(k-1) &\text{otherwise}
\end{cases}}
Conditions \eqn{75} and \eqn{76} are satisfied with $\domain =
\badchoice$ and $\atom = \avoid$. Also, for all $(y,b) \in \badchoice$,
\display{79}{\avoid(y,b) = f \implies
 \exists m \, [\, f(m) = y \mand \forall k \! \ge
    \! m \; \zeta(f(k)) \neq b \, ]}
}

\begin{proof}
First it should be clarified how \eqn{78} uniquely defines the function $f =
\avoid(x,a)$, where $\langle z_0, \dotsc, z_n \rangle$ be the path
from $\nothing$ to $x$. For $k < n$, the first clause of \eqn{78} is an
explicit definition. For $k > n$, \eqn{78} is equivalent to the following
definition by recursion on $\Nat$, having $i$ as the variable on which
recursion is performed.
\display{80}{\begin{aligned} f(n+0) &= z_n\\ f(n+i+1) &= \begin{cases}
      \min \left\{ y \mid y \in Y(f(n+i))
      \right.&\\ \strut \hspace{4em}\left. \mand \zeta(y) \neq a
      \right\} &\text{if\ } Y(f(n+i)) \neq \emptyset\\ f(n+i)
      &\text{otherwise}
\end{cases}
\end{aligned}}

To verify that condition \eqn{75} holds, suppose that, for some $x$
and $a$, that $n$ is the depth of $x$ and $(x,a) \in \badchoice$ and
$\avoid(x,a) = f$ and the antecedent of \eqn{75} is satisfied at
$(x,a)$. That is, suppose that, for some $y$, $y \sms x$. Since $(\X,
\wms)$ is a tree, $\{ w \mid y \wms w \sms x \}$ is finite, so it is
well ordered by $\wms^{\mathord{-}1}$, and therefore has a least
specific element. That is, $Y(x) \neq \emptyset$. Since $\zeta$
is ID consistent and $(x,a) \in \badchoice$, there is some $z \in
Y(x)$ such that $\zeta(z) \neq a$. Thus $f(n+1)$ is defined by the
first alternative of \eqn{80}, that $f(n+1) = \min \{ y \mid y \in
Y(f(n+i)) \mand \zeta(y) \neq a \}$. Assume without loss of generality
that $z = \min \{ y \mid y \in Y(f(n+i)) \mand \zeta(y) \neq a \}$, so
that $f(n+1) = z$. Then $(z,a) \in \badchoice$. Let $\avoid(z,a) =
g$. Then $f(k) = g(k)$ for $k \le n$, and, by induction on $\Nat$
using \eqn{80}, $\forall i\! \in \! \Nat \; f(n+i) = g(n+i)$. That is,
$\avoid(x,a) = \avoid(z,a)$, establishing that the consequent of
\eqn{45} holds for $x$ and $a$.

To verify that condition \eqn{76} holds, suppose that $(x,a) \in
\badchoice$ and $(y,b) \in \badchoice$ and, without loss of
generality, that the depths of $x$ and $y$ are $i$ and $j$
respectively and that $i \le j$. Suppose also the negation of the
consequent of \eqn{76} so that, for some $f$, $f \in \atom(\domain
\cap (\{ x \} \times A)) \cap \atom(\domain \cap (\{ y \} \times
A))$. By \eqn{78}, $f(i) = x$ and $f(j) = y$.  By induction on $n-m$,
if $m \le n$ and $\avoid(x,a) = f$, then $f(n) \wms f(m)$.  Therefore
$y \wms x$, so not $x \incompat y$. Condition \eqn{76} holds by
contraposition.

Condition \eqn{79} follows immediately from \eqn{78} by using induction on
$k$ to prove $\avoid(y,b) = f \implies [f(m) = y \mand \forall k \!
  \ge \! m \; \zeta(f(k)) \neq b]$, taking $m$ to be the depth of $y$.
\end{proof}

\corol S{Define $\badchoice$ by \eqn{52} and define $\avoid \from
  \badchoice \to \X^\Nat$ by setting $\avoid(x,a) = f$, where $f$
  satisfies \eqn{78}. Define $\range = \avoid(\badchoice)$. Define
  $\sample[\zeta]$, the graph of $\avoid$, by \eqn{59}. Define $\br[\zeta]
  \from \X \to 2^\range$ by \eqn{77}, with $\atom = \avoid$ and
  $\domain = \sample[\zeta]$. Define $\event[\zeta] \from \X \to
  2^{\sample[\zeta]}$ by
\display{81}{\event[\zeta](x) = \{ (z,a,f) \mid (z,a,f) \in
  \sample[\zeta] \mand \exists n \; f(n) = x \}}
Then $\event[\zeta]$ is an embedding of $(\X, \wms)$ in
$(\sample[\zeta], 2^{\sample[\zeta]})$.}

\begin{proof}
As $\domain$ and $\range$ are defined in the corollary, condition
\eqn{74} is clearly satisfied. By \rslt A, conditions \eqn{75} and
\eqn{76} are also satisfied. Therefore, by \rslt L, $\br[\zeta]$ is an
embedding. Finally, by \rslt M and \eqn{81}, $\event[\zeta]$ is also
an embedding.
\end{proof}

With the foregoing preparation, the proof of \rslt Q can now be
completed.

\begin{proof}[Proof of \Rslt Q]
As in \rslt S, define $\badchoice$ by \eqn{52}, $\avoid$ by setting
$\avoid(x,a) = f$ where $f$ satisfies \eqn{52}, $\range =
\avoid(\badchoice)$, $\sample[\zeta]$ by \eqn{59}, $\br[\zeta]$ by
\eqn{77} with $\atom = \avoid$ and $\domain = \sample[\zeta]$, and
$\event[\zeta]$ by \eqn{81}.

Define state-contingent utility functions, $\{ u_a \from
\sample[\zeta] \to \Re \}_{a \in A}$, by
\display{82}{u_a(x,b,f) = \begin{cases} 1, & \text{if\ } \exists i \;
    \exists j \; [i \le j \mand f(i) = x \mand \mand \zeta(f(j)) =
      a]\\ 0, & \text{otherwise}
\end{cases}}

Define subsets $S$ and $D$ of $\sample[\zeta]$ by
\display{83}{\begin{aligned}
(y,b,f) \in S \iff [(y,b,f) &\neq \omega \mand\\
&\exists i \; \exists j \; [i \le j \mand y = f(i) \mand x = f(j)]]
\end{aligned}}
and
\display{84}{\begin{aligned}
(y,b,f) \in D \iff [y &\neq x \mand\\
&\exists i \; \exists j \; [i < j \mand x = f(i) \mand y = f(j)]]
\end{aligned}}

From this point forward, \rslt Q is proved by repeating word-for-word
the proof of \rslt H, except for the handful of amendments that will
be specified below. The overall principle underlying these amendments
is that, wherever the proof of \rslt H envisions that the range of
$\avoid$ is $\leaves$, the proof of \rslt Q envisions rather that the
range of $\avoid$ is a subset of $\X^\Nat$. In particular, where a
statement that $z \wms w$ (where $z = \avoid(y,b)$ is a leaf) appears
in the proof of \rslt H, that statement should be replaced in the
proof of \rslt Q by a statement that $\exists i \; f(i) = w$, where $f
= \avoid(y,b)$.
\begin{itemize}
\item 
The assertion that $\sample[\zeta]$ is finite, in the sentence
immediately preceding the equations \eqn{60}, should be replaced by an
assertion (justified by \rslt{R}) that $\sample[\zeta]$ is
countable.\footnote{Since $\X$ is countable and $A$ is finite,
  $\badchoice$ is countable; and since there is a 1--1 correspondence
  between $\badchoice$ and $\sample[\zeta]$, $\sample[\zeta]$ is
  countable.} Since every function of a countable domain is measurable
with respect to the power set of the domain, the utility functions
$u_a$ defined in \eqn{82} are measurable functions.
\item
Definitions \eqn{61} and \eqn{62} are replaced by definitions \eqn{82}
and \eqn{81} respectively. Definition \eqn{65} is replaced by
definitions \eqn{83} and \eqn{84}.  Subsequent references to \eqn{61},
\eqn{62}, and \eqn{65} must also be revised.
\item
A reference to condition \eqn{54}, at the beginning of the paragraph
where condition \eqn{67} is stated, should be replaced by a reference
to condition \eqn{79}. Following the reference to \eqn{54} is a
paraphrase, `$\avoid(x,a) \wms x\/$', of a clause of the
condition. That paraphrase should be amended to `$\exists i \, f(i) =
x$, where $f = \avoid(x,a)\/$'.
\item
At the beginning of the paragraph where inequality \eqn{68} is stated,
`$\psi = (y,b,z)\/$' should be amended to `$\psi = (y,b,f)\/$'.  In a
paraphrase of \eqn{65} in the sentence immediately preceding the
statement of condition \eqn{68}, the clause stating that `$z \wms x
\wms y\/$', should be amended to `$\exists i \; \exists j \; [i \le j
  \mand y = f(i) \mand x = f(j)]\/$'.
\end{itemize}

\end{proof}

\Section c{Concluding remarks}

\Subsection g{Antecedent theories of subjective conditional probability}

While Savage does not formalize the notion of an agent making a choice
after having received evidence, he regards his theory as ``giv[ing] an
atemporal analysis of [a] temporally described decision situation: the
person must decide between [two alternatives] after he
\dots\ observes, whether [an event] obtains[.]''
\citeyearpar[p.\ 23]{Savage-1972} Savage remarks in
\citep{AumannSavage-1990} that his discussion of \emph{small
  worlds} in \citeyearpar[section 5.5]{Savage-1972} reflects his
awareness of difficulties implicit in that interpretation,
particularly regarding the specification that consequences have
utilities that are independent of the state of the world.
\citet{Dreze-1961} and Aumann, in \citep{AumannSavage-1990}, fully
articulated those difficulties, and noted that the
uniqueness of probabilities and (up to a positive affine
transformation) utilities derived by Savage would not survive deleting
the state-independent-utility assumption. Slightly prior to Aumann's
critique, several researchers had provided
models of, and representation theorems for, explicit conditioning of
subjective expected utility on events. (\citet{KrantzEtAl-1971}
exposit their own and other researchers' contributions.) These models,
and also a closely related model due to \citet{Skiadas-1997b}, take
acts and states of the world to be primitive, and treat the
consequence of taking an act in a state of the world to be simply the
ordered pair of the act and the state; while an alternate model due to
\citet{Karni-2015} takes acts and consequences to be primitive, and
constructs states of the world as functions from acts to
consequences. In all cases, and also in the models of dynamic
consistency and consequentialism mentioned in footnote
\ref{fn:dynamic}, evidential states are identified with events, which
are modeled as sets of states of the world, and those sets are assumed
to form a Boolean algebra.

\Subsection d{Finitely versus countably additive probabilities}

As \citet{SchervishEtAl-1984} document, Bayesians have been divided on
whether finitely or countably additive probability measures provide
the more conceptually appropriate representation of subjective
probabilities. While \sceup ity has been defined here in terms of
countably additive probabilities, all of the results established in
this article are robust to changing to a definition in terms of
finitely additive probabilities. The choice between the two frameworks
has three implications. First, finitely additive measures can be
defined on all subsets of a set, while, for uncountable sets, the
domain of a countably additive measure must generally be a strict
subset of the power set of that set. All of the examples formulated in
this article, as well as the constructions utilized in the proofs of
\prp{H} and \prp{Q}, involve measures that are defined on the full
power set of a countable set, so that distinction is moot. Second,
integration with respect to a countably additive measure is well
defined for a class of unbounded functions, while integration with
respect to a finitely additive measure is defined for bounded
functions only. Again, boundedness of utility functions has been
imposed in the definition of \sceup ity for reasons independent of
finite additivity, the distinction is moot. Third, even in cases where
a countably additive measure can be defined on the full power set of
an infinite set, there are finitely additive measures that are not
countably additive. Since the examples and propositions of this
article concern the existence of probability measures rather than the
properties shared by all probability measures, and since countably
additive measures have been constructed or proved to exist and those
measures are \emph{a fortiori\/} finitely additive, the distinction is
moot once again.

\Subsection e{The parsimonious theory is weak}

One way of thinking about the strength of theories is that a theory is
strong if it has few models, and especially if all its models resemble
the intended model in salient respects. In \sec{w}, an evidential
structure was introduced to represent the infinitely-repeated tossing
of a coin, and the intended interpretation of that structure was
described. Specifically, the e-structure is intended to describe an
uncountable probability space with a non-atomic probability
measure. It is obvious how to specify a plan to represent the choices
of a risk-neutral agent who has an opportunity to bet on whether or
not the coin is biased towards `heads' (that is, on whether or not the
long-run average of the fraction of tosses that land `heads' will
exceed $1/2$): bet on heads in an evidential state (that is, after
having observed a finite number of tosses) if, and only if, more than
half of the tosses have landed `heads'. There is an intuitive class of
intended models that show this plan to be \sceup. That is, two coins,
of bias $(1+x)/2$ and $(1-x)/2$ (for $0 < x < 1$) are placed in an
urn, and then one of them is drawn at random (with probability $1/2$,
according to the agent's beliefs) and is tossed repeatedly. Like the
intended model of a sequence of tosses of a known coin, the models in
this class involve non-atomic measures on uncountable sets. However,
by \rslt Q, the plan is also rationalized by a model involving a
purely atomic probability measure on a countable set. This phenomenon,
that a theory with an uncountable intended model also has an countable
model, is reminiscent of the L\"{o}wenheim-Skolem theorem
\citep[theorem 3.5]{BellSlomson-2006} regarding first-order
theories. The present theory of \sceup ity is a second-order theory,
since the definition of \sceup ity has the form, ``There exist
functions, $u_a$, such that\dots''. Nonetheless, it is a weak theory
in precisely the respect that the L\"{o}wenheim-Skolem theorem is
typically understood to exhibit the weakness of first-order theories.

\Subsection h{Embedding an infinite e-structure in a Boolean algebra}

\Rslt{S} is an analogue of \rslt{C} for countable e-structures that
are trees. In fact, every e-structure, even an uncountable
e-structure that are not a tree, can be embedded in a Boolean
algebra. This fact follows from theorems 7.11 (Stone's representation
theorem) and theorem 14.10 of \citet{Jech-2002}, and was
used by \citet{bg-2016} to prove \rslt{Q}. It might be useful for
generalizing the present theory.

\Subsection f{\uppercase\sceup ity corresponds to na\"{i}ve Bayesian rationality}

In the Bayesian subjectivist theory of inference (sketched by Savage
in \citeyearpar[chaptes 6,7]{Savage-1972}), an evidential state
corresponds to an \emph{observation.} The agent observes one thing or
another, so distinct evidential states are reached in disjoint
events. Therefore (assuming that something or other will be observed),
the evidential states correspond to a partition of the states of the
world. Thus their probabilities sum to one and, in particular, the sum
of probabilities of several evidential states cannot exceed one. In the theory
developed in the present article, the event associated with an
evidential state, $x$, is $\event(x)$. In general, the events
associated with distinct evidential states are not
disjoint. Consequently the sum of their probabilities can exceed
one. This possibility plays a role in the proof of \rslt{s}, for
example, where it is specified that $\Pr(\event(x_1)) =
\Pr(\event(x_2)) = 5/7$. That is, the subjective probabilities
ascribed to the \sceu agent in the example violate a constraint that a
Bayesian analysis would impose on the agent's beliefs.

The questions then arise, under what conditions on an evidential
structure do the probabilities ascribed to the agent by every
embedding of the e-structure into a probability space necessarily
satisfy the Bayesian constraints, and under what conditions does there
at least exist some embedding for which the ascribed probabilities
satisfy the Bayesian constraints. \citet{bg-2018} study these
questions. They find that only trivial e-structures, in which every
e-state in $\X \setminus \{ \nothing \}$ is immediately more specific
than $\nothing$, necessarily induce probabilities that reflect an
underlying Bayesian model. They also find that, for a class of
e-structures that strictly includes trees, some embedding into a
probability space induces probabilities that reflect an underlying
Bayesian model. Those findings provide an argument from a Bayesian
perspective, complementary to the present results, for the suggestion,
offered at the end of \sec{1}, that the version of \sceup ity
formulated in this article might better be regarded as a heuristic
criterion than as an intrinsic definition of rational evidence-based
choice.

\bibliographystyle{plainnat}

\strut \end{document}